\documentclass[12pt]{amsart}
\usepackage{amssymb, latexsym, amsthm}

\newtheorem{thm}{Theorem}[section]

\newtheorem{cor}[thm]{Corollary}
\newtheorem{defn}[thm]{Definition}
\newtheorem{exmpl}[thm]{Example}
\newtheorem{lem}[thm]{Lemma}
\newtheorem{prop}[thm]{Proposition}
\newtheorem{rem}[thm]{Remark}

\newtheorem{algo}[thm]{Algorithm}

\def\ie{{i.e.}}
\def\eg{{e.g.}}
\def\cf{{\it cf. \/}}

\def\lam{{\lambda}}
\def\uni{{\varpi}}

\def\color{\varrho}
\def\n{{\nu}}
\def\s{\sigma}

\def\A{{\mathbb A}}
\def\C{{\mathbb C}}
\def\F{{\mathbb F}}

\def\N{{\mathbb N}}

\def\Z{{\mathbb Z}}

\def\AA{{\mathcal A}}
\def\B{{\mathcal B}}

\def\O{{\mathcal{O}}}

\def\Vals{{\mathcal{V}}}

\def\tG{{\tilde{G}}}

\DeclareMathAlphabet{\mathscr}{OT1}{pzc}{m}{it}

\newcommand\tensor[1][{}]{{\otimes_{#1}}}

\def\ra{{\rightarrow}}
\def\lra{{\longrightarrow}}

\def\hra{{\,\hookrightarrow\,}}
\def\minusset{{-}}

\def\sub{\subseteq}

\def\({\left(}
\def\){\right)}
\def\isom{{\;\cong\;}}
\def\normal{{\unlhd}}
\def\semidirect{{\ltimes}}
\def\normali{{\lhd}}
\def\co{{\,{:}\,}}
\def\divides{{\,|\,}}
\newcommand\suchthat{{\,:\ \,}}
\newcommand\subjectto{{\,|\ }}

\newcommand\comp[1]{{{#1}^{\operatorname{c}}}}

\DeclareMathOperator{\Ker}{Ker}

\DeclareMathOperator{\diag}{diag}

\newcommand\Cayley[2]{{\operatorname{Cay}(#1;#2)}}

\DeclareMathOperator{\End}{End}

\DeclareMathOperator{\mychar}{char} \DeclareMathOperator{\Cent}{Z}
\DeclareMathOperator{\tr}{tr}
\newcommand{\Norm}[1][]{{\operatorname{N}_{#1}}}
\DeclareMathOperator{\Gal}{Gal}

\DeclareMathOperator{\Br}{Br}

\newcommand\cond[2][!]{{\operatorname{cond}_{\if!#1\relax\else{\comp{#1}}\fi}(#2)}}

\renewcommand\L[2]{{\operatorname{L}^{#1}(#2)}}
\newcommand\M[1][d]{{\operatorname{M}_{#1}}}
\newcommand\GL[1][d]{{\operatorname{GL}_{#1}}}
\newcommand\PGL[1][d]{{\operatorname{PGL}_{#1}}}

\newcommand\SL[1][d]{{\operatorname{SL}_{#1}}}

\newcommand\PSL[1][d]{{\operatorname{PSL}_{#1}}}

\newcommand\abs[2][F]{|{#2}|_{#1}}

\newcommand{\set}[1]{{\{#1\}}}
\newcommand{\card}[1]{{\left|{#1}\right|}}
\newcommand\ideal[1]{{\left<{#1}\right>}}
\newcommand\sg[1]{{\ideal{#1}}}

\newcommand\dimcol[2]{{[{#1}\!:\!{#2}]}}

\newcommand\db[1]{{(\:\!\!({#1})\:\!\!)}}
\newcommand\mul[1]{{#1^{\times}}}
\newcommand\md[2][d]{{\mul{#2}/{\mul{#2}}^{#1}}}

\newcommand\He[2][G]{{\operatorname{H}(#1,#2)}}

\newcommand\valF{{\nu_0}}
\newcommand\Fy[1][q]{{\F_{#1}\db{y}}}

\newcommand\Tref[1]{{Theorem \ref{#1}}}

\newcommand\Pref[1]{{Proposition \ref{#1}}}

\newcommand\Cref[1]{{Corollary \ref{#1}}}
\newcommand\Dref[1]{{Definition \ref{#1}}}
\newcommand\Rref[1]{{Remark \ref{#1}}}

\newcommand\eq[1]{{(\ref{#1})}}

\newcommand\Eq[1]{{Equation \eq{#1}}}

\newcommand\defin[1]{{\it{#1}}}
\long\def\half#1\halved{{\footnotesize{#1}}}
\long\def\forget#1\forgotten{}

\newcommand\dc[3]{{{#1}\backslash{#2}/{#3}}}
\newcommand\dom[2]{{{#1}\backslash{#2}}}

\newcommand\paper[6]{{{#1},\ {\it{#2}},\ {#3}\ {#4},\ {#5},\ ({#6}).}}
\newcommand\book[4]{{{#1},\ {{#2}},\ {#3},\ {#4}.}}

\newcommand\dd[4][!]{\abs[#2]{\:\!\det{\if!#1\relax\:\!\!\else(#1)\fi}}^{#3} #4}

\newcommand\yy[1]{{{#1}_{{\:\!} 1 {\:\!\!} / {\;\!\!} y}}}

\newcommand\binomq[3]{{\genfrac{[}{]}{0pt}{2}{#1}{#2}_{#3}}}
\newcommand\con[3]{{{#1}({#2},{#3})}}
\def\Aspec{{\mathfrak S}}

\def\pp{p}
\def\chix{x}

\newif\ifXY
\XYtrue

\ifXY
\usepackage{xy}
\fi \ifXY \xyoption{all} \fi

\begin{document}

\title
[Explicit constructions of Ramanujan Complexes] {Explicit
constructions of Ramanujan Complexes}

\def\HUJI{Inst. of Math., Hebrew Univ., Givat Ram, Jerusalem 91904,
Israel}
\def\BIU{Dept. of Math., Bar-Ilan University, Ramat-Gan 52900,
Israel}
\def\YALE{Dept. of Math., Yale University, 10 Hillhouse Av., New-Haven CT
06520, USA}

\author{Alexander Lubotzky}
\address{\HUJI}
\email{alexlub@math.huji.ac.il}

\author{Beth Samuels}
\address{\YALE
}
\email{beth.samuels@yale.edu}

\author{Uzi Vishne }
\address{\YALE}
\curraddr{\BIU}

\email{vishne@math.biu.ac.il}

\thanks{This research was supported by the NSF and the BSF (U.S.-Israel)}

\renewcommand{\subjclassname}{
      \textup{2000} Mathematics Subject Classification}
\date{Received: Aug. 31, 2003, Revised: Feb. 25, 2004}

\begin{abstract}
In this paper we present for every $d \geq 2$ and every local
field $F$ of positive characteristic, explicit constructions of
Ramanujan complexes which are quotients of the Bruhat-Tits
building $\B_d(F)$ associated with $\PGL[d](F)$.
\end{abstract}

\maketitle

\section{Introduction}

\def\ww{{y}}

In \cite{paperI} we defined and proved the existence of Ramanujan
complexes, see also \cite{Ballantine}, \cite{CSZ} and \cite{Li2}.
The goal of this paper is to present an explicit construction of
such complexes.

Our work is based on the lattice constructed by Cartwright and
Steger \cite{CS1}. This remarkable discrete subgroup $\Gamma$ of
$\PGL[d](F)$, when $F$ is a local field of positive
characteristic, acts simply transitively on the vertices of the
Bruhat-Tits building $\B_d(F)$, associated with $\PGL[d](F)$. By
choosing suitable congruence subgroups of $\Gamma$, we are able to
present the $1$-skeleton of the corresponding finite quotients of
$\B_d(F)$ as Cayley graphs of explicit finite groups, with
specific sets of generators. The simplicial complex structure is
then defined by means of these generators.

Let $\binomq{d}{k}{q}$ denote the number of subspaces of dimension $k$
of $\F_q^d$.
\begin{thm}\label{mainIntro}
Let $q$ be a prime power, $d \geq 2$, $e \geq 1$ ($e > 1$ if $q =
2$).

Then, the group $G = \PGL[d](\F_{q^e})$ has an (explicit) set $S$
of $\binomq{d}{1}{q} + \binomq{d}{2}{q} + \cdots +
\binomq{d}{d-1}{q}$ generators, such that the Cayley complex of
$G$ with respect to $S$ is a Ramanujan complex, covered by
$\B_d(F)$, when $F = \F_q\db{y}$.
\end{thm}
The Cayley complex of $G$ with respect to a set of generators $S$
is the simplicial complex whose $1$-skeleton is the Cayley graph
$\Cayley{G}{S}$, where a subset of $i+1$ vertices is an $i$-cell
iff every two vertices comprise an edge.
The generators in \Tref{mainIntro} are explicitly given in Section
\ref{sec:explicit}.

In the case $d = 2$ there are two types of Ramanujan graphs,
bi-partite, and non-bi-partite (\cite{LPS}, \cite{Morg},
\cite{alexbook}). Here too, given $r$ dividing $d$, we produce
$r$-partite complexes, by proving an analog of Theorem
\ref{mainIntro} for every subgroup of $\PGL[d](\F_{q^e})$
containing $\PSL[d](\F_{q^e})$.
There is also a version of the theorem for subgroups of
$\PGL[d](L)$ where $L$ is a finite local ring.

The paper is organized as follows: In Section \ref{sec:buildings}
we describe affine buildings of type $\tilde{A}_{d}$ in general,
in the language of
$\O$-sublattices of $F^{d}$ (where $\O$ is the valuation ring of
$F$). The Cartwright-Steger group $\Gamma$ is constructed in
Sections \ref{sec2}--\ref{sec:gamma}. Our construction slightly
differs from theirs, but is more convenient for the computations
to follow. The simply transitive action of $\Gamma$ on the
building is used in Section \ref{geom} to describe the defining
relations of $\Gamma$.

In Section \ref{sec:finq} we present and analyze the finite
Ramanujan quotients. In particular, we use \cite{paperI} to show
that the complexes constructed here are Ramanujan. We should note
that the proof of this result in \cite{paperI} relies on what is
called the global `Jacquet-Langlands correspondence' for function
fields (a correspondence between automorphic representations of a
division algebra and of $\GL[d]$, \cf\ \cite[Thm~VI.1.1]{HT} for
the characteristic zero case). This correspondence in the function
field case is also considered to be true by experts, as the main
ingredients of the proof are known; though the task of writing
down a complete proof has not been carried out yet.

Originally, $\Gamma$ is a group of $d^2\times d^2$ matrices over a
ring $R$. Section \ref{splitx} provides an explicit embedding of
$\Gamma$ into $d\times d$ matrices over a finite extension of $R$.
In Section \ref{sec:explicit} this embedding is refined to
identify finite quotients of $\Gamma$ with subgroups of
$\PGL[d](L)$ which contain $\PSL[d](L)$, where $L$ is a finite
local ring; a detailed algorithm is given, with an example in
Section \ref{sec:example}. In particular the generators of
$\Gamma$ are given as $d\times d$ matrices over $\F_q[x]$. For the
convenience of the reader, we include a glossary in the final
section.

{\it Added in Proof:}
We recently learned that Alireza Sarveniazi \cite{Ali}
has also given an explicit construction of Ramanujan complexes.

\section{Buildings}\label{sec:buildings}

To every reductive algebraic group over a local field one can
associate a building, which is a certain simplicial complex, on
which the group acts (see \cite{Ronan}).
This complex plays the role of a symmetric space for Lie groups.

Recall that a complex is a structure composed of $i$-cells, where
the $0$-cells are called the vertices, and an $i$-cell is a set of
$i+1$ vertices. A complex is simplicial if every subset of a cell
is also a cell. The $i$-skeleton is the set of all $i$-cells in
the complex. Buildings are in fact clique complexes, which means
that a set of $i+1$ vertices is a cell iff every two vertices form
a $1$-cell. This property holds for quotient complexes, which will
be the subject of Section \ref{sec:finq}.

We will now describe the affine building associated to $\PGL(F)$,
where $F$ is a local field.
These are called `buildings of type $\tilde{A}_{d-1}$' because of
the Dynkin diagram of the associated Weil group (which is
isomorphic to $S_d \semidirect \Z^{d-1}$). Let $\O$ denote the
valuation ring of $F$; choose a uniformizer $\uni$ (so for $F =
\Fy$, $\O = \F_q[[y]]$ and $\uni = y$), and assume $\O/\uni \O =
\F_q$. Consider the $\O$-lattices of full rank in $F^d$. For every
lattice $L$, $\uni L$ is a sublattice, and as $\F_q$-vector spaces
$L/\uni L \isom \F_q^d$. We define an equivalence relation by
setting $L \sim s L$ for every $s \in \mul{F}$.
\begin{rem}\label{inver}
Let $\sg{\uni}$ denote the multiplicative subgroup of $\mul{F}$
generated by $\uni$, and let $\mul{\O}$ be the invertible elements
of $\O$. Then $\mul{F} = \sg{\uni}\cdot \mul{\O}$.
\end{rem}
Since $L = s L$ for any element $s \in \mul{\O}$, the equivalence
classes have the form $[L] = \set{\uni^{i}L}_{i\in \Z}$. Let
$\B^0$ be the graph whose vertices are the equivalence classes.
There is an edge from $[L]$ to a class $\chix \in \B^0$, iff there
is a representative $L' \in \chix$ such that  $\uni L \subset L'
\subset L$. Notice that this is a symmetric relation, since then
$\uni L' \subset \uni L \subset L'$.

The vertices of $\B^0$ form the $0$-skeleton of a complex $\B$,
and the edges are the $1$-skeleton, $\B^1$. As $i$-cells of $\B$
we take the complete subgraphs of size $i+1$ of $\B^0$, which
correspond to flags
$$\uni L_0 \subset L_i \subset \dots \subset L_1 \subset L_0;$$
the $i$-skeleton is denoted $\B^{i}$. It
immediately follows that $\B$ has $(d-1)$-cells (corresponding to
maximal flags in quotients $L/\uni L$). It also follows that there
are no higher dimensional cells.

The group $\GL(F)$ acts transitively on  lattices by its action on
bases. moreover note that the action preserves inclusion of
lattices. We call $L_0 = \O^d \sub F^d$ the
\defin{standard lattice}
.
If $\tau \in \GL(F)$ has entries in $\O$, then $\tau L_0 \sub
L_0$. The stabilizer of $L_0$ in $G$ is thus the maximal compact
subgroup $\GL(\O)$. According to the definition of the equivalence
relation, the scalar matrices of $\GL(F)$ act trivially on $\B$,
so the action of $\GL(F)$ induces a well defined action of
$\PGL(F)$ on (the vertices of) $\B$, which is easily seen to be an
action of an automorphism group. Again, the stabilizer of $[L_0]$
is the maximal compact subgroup $\PGL(\O)$. The set of vertices
can thus be identified with $\PGL(F)/\PGL(\O)$.

Since the only ideals of $\O$ are powers of $\ideal{\uni}$,  The
Invariant Factor Theorem for $F$ asserts that any matrix in
$\GL[d](F)$ can be decomposed as $a g a'$ for $a,a' \in
\GL[d](\O)$ and $g = \diag(\uni^{i_1},\dots,\uni^{i_d})$, where
$i_1 \leq \dots \leq i_d$ are integers. If $L = a g a' L_0 = a g
L_0$, then $\uni^{-i_1} L = a (\uni^{-i_1} g)L_0 \sub L_0$, and on
the other hand $\uni^{-i_d} L = a (\uni^{-i_d} g)L_0 \supseteq
L_0$, so any two lattices of maximal rank are commensurable.
Moreover in this case, $L_0/\uni^{-i_1}L$ is annihilated by
$\uni^{i_d-i_1}$, and so is a module over $\O / \uni^{i_d-i_1}
\O$, a local ring of order $q^{i_d-i_1}$. In particular
$\dimcol{L_0}{\uni^{-i_1} L}$ is a $q$ power.

This basic fact allows us to define a color function $\color \co
\B^0 \ra \Z/d$ by $\color(L) = \log_q \dimcol{L_0}{\uni^iL}$ for
large enough $i$; this function is well defined since
$\dimcol{\uni^i L}{ \uni^{i+1} L} = q^d$. Also notice that
$\color([\tau L_0]) = \valF(\det(\tau)) \pmod{d}$, when $\valF$ is
the valuation of $F$. This shows that $\SL[d](F)$ is color
preserving, while on the other hand, $\tau_\uni =
\diag(\uni,1,\dots,1)$ has determinant $\uni$, so
$\color(\tau_\uni(L)) = \color(L) + 1$ for every $L$. It follows
that $\GL(F)$ acts transitively on colors.

The color is additive in the sense that if $L'' \sub L'$, then
$\color(L'') = \color(L')+\log_q\dimcol{L'}{L''}$. Similarly if $L
\sub L_0$ and $\tau \in \GL(F)$, then $\color(\tau L) =
\color(\tau L_0) + \color(L)$.

The colors provide us with $d$ Hecke operators, defined on
functions of $\B^0$ by summation over the neighbors of fixed
color-shift:
$$A_kf(x) = \sum_{y\sim x,\, \color(y)-\color(x)\equiv k} f(y).$$
These operators generate the Hecke
algebra $\He[\PGL(F)]{\PGL(\O)}$ (see \cite[Sec.~2]{paperI} for
more details).

For $1 \leq t < d$, let $\B[t]$ denote the graph defined on the
vertices $\B^0$, with the edges $(\chix,\chix') \in \B^1$ for
which there are $L \in \chix$ and $L' \in \chix'$ such that $L'
\sub L$ and $\dimcol{L}{L'} = q^t$ (in particular,
$\color(\chix)-\color(\chix') = t$).
\begin{rem}\label{close}
If $L'\subset L$ is a sublattice of index $q$, then $[L],[L']$ are
connected in $\B^1$.
\end{rem}
\begin{proof}
We need to prove that $\uni L \sub L'$, but this is obvious since
$L/L'$ is annihilated (as an $\O$-module) by multiplication by
$\uni$.
\end{proof}

For every $L'\sub L$ there is a composition series of sublattices
$L'= L^{(m)} \sub L^{(m-1)} \sub \cdots \sub L^{(0)} = L$, such
that $\dimcol{L^{(i)}}{L^{(i+1)}} = q$. It follows that $\B[1]$ is
a connected (directed) subgraph of $\B$. In fact, $\B[1]$
determines $\B^1$, and thus all of $\B$:
\begin{prop}\label{partofd}
Vertices $\chix,\chix'$ are connected in $\B^1$ iff there is a
chain $\chix_0,\chix_1,\dots,\chix_d \in \B^0$ such that $\chix_0
= \chix_d = \chix$, $(\chix_i,\chix_{i+1}) \in \B[1]$ for $i =
0,\dots,d-1$, and such that $\chix' \in
\set{\chix_1,\dots,\chix_{d-1}}$.
\end{prop}
\begin{proof}
First assume that such a chain exists, and choose representatives
$L_i \in \chix_i$ with $L_d = \uni L_0$. By the definition of
$\B[1]$, we may assume that $\dimcol{L_{i}}{L_{i+1}} = q$ and then
$L_{d} \subset \cdots \subset L_2 \subset L_1 \subset L_0$. In particular
$([L_i],[L_0]) \in \B^1$ for every $i$.

On the other hand, if
\begin{equation}\label{LLL}
 \uni L_0 \subset L' \subset L_0
\end{equation}
are lattices, we can lift a maximal flag in
$L/\uni L_0 \isom \F_q^d$ to a maximal chain of lattices refining
\eq{LLL}, resulting in a chain $\chix_0,\dots,\chix_d$.
\end{proof}

\begin{cor}
If $(\chix,\chix') \in \B[t]$, then there is a path of length $t$
in $\B[1]$ from $\chix$ to $\chix'$.
\end{cor}

{}Using this criterion, it is easy to see that if the greatest
common divisor $(t',d)$ equals $t$, and $(\chix,\chix') \in
\B[t']$, then there is a path from $\chix$ to $\chix'$ in $\B[t]$.
In particular, $\B[t']$ has the same connected components as
$\B[t]$. The final result of this section is not needed in the
rest of the paper.

We thank the referee for some simplification in the proof of the
following proposition.

\begin{prop}\label{Bt}
Let $t$ be a divisor of $d$. If $t$ divides
$\color(\chix')-\color(\chix)$, then there is a path from $\chix$
to $\chix'$ in the (directed) graph $\B[t]$ defined above.
\end{prop}
\begin{proof}
We may assume $\chix = [L_0]$, and $\chix' = [L]$ with $L = aga'L_0 = agL_0$,
$a,a' \in \GL[d](\O)$ and $g =
\diag(\uni^{i_1},\dots,\uni^{i_d})$. We may assume $0 \leq i_1 \leq \dots \leq
i_d$, and in particular $L \sub L_0$. If the claim is true for $g
L_0 \sub L_0$, then acting with $a$ we get a chain from $a L_0
=L_0$ to $L$; we can thus assume $a = 1$.

Before going on by induction, we multiply $L$ by a power of
$\uni$ so that
\begin{equation}\label{strange}
i_1 \geq \frac{t}{d-t}(i_d-i_1).
\end{equation}
Now, by
assumption, $i_1+\cdots + i_d = rt$ for some $r \in \N$.

{\it Case 1.} If $i_{d-t+1} > i_1$, then we can lower $t$ of the entries
$i_2,\dots,i_d$, each by $1$, keeping the increasing
order. Let $i_1',\dots,i_d'$ denote the resulting values and $g' =
\diag(\uni^{i_1'},\dots,\uni^{i_d'})$.
Then $g L_0 \subset g' L_0 \sub L_0$ with
$\dimcol{g' L_0}{g L_0} = q^t$ and $\uni g' L_0 \sub g L_0$, so
$(g'L_0,gL_0) \in \B[t]$. Since the condition $i_1 \geq
\frac{d}{t}(i_d-i_1)$ still holds (as $i_1$ was not changed),  we
are done by induction on $r$.

{\it Case 2.} Now assume $i_{d-t+1} = i_1$; let $j$ be maximal with $i_j =
i_1$. If $j = d$ then $i_1 = \dots = i_d$ so $L \equiv L_0$ and we are done.
Therefore assume $j < d$. Of course $j \geq d-t+1$.
If $i_1 = 0$ then $i_d = 0$ too by the assumption \eq{strange}, so $L = L_0$
and again we are done. Assume $i_1 > 0$ and $j < d$.
In the first step we can lower $i_{j+1},\dots,i_d$, but we also have to lower
$i_1,\dots,i_{t-(d-j)}$ in order to change exactly $t$ components. We continue lowering the highest entries, using the remaining $d-t$ entries
$i_{t-(d-j)+1},\dots,i_{j}$ whenever necessary. In each of the
next $d/t-1$ steps, $i_d$ is lowered by one with $i_1$
unharmed---so at the end our condition \eq{strange} is still met. The
$d/t$ modules we constructed form a chain in $L_0$, climbing from $L$
using $\B[t]$ steps, and ending with indices which now sum up to $rt-d$,
so again we are done by induction.
\end{proof}

Vertices $\chix$ and $\chix'$ for which $(t,d)$ does not divide
$\color(\chix')-\color(\chix)$ cannot be connected in $\B[t]$, so
we proved that $\B[t]$ has $(t,d)$ connected components for every
$1\leq t < d$.

\section{The arithmetic lattice}\label{sec2}

Let $\F_q$ denote the field of order $q$ (a prime power), and
$\F_{q^d}$ the extension of dimension $d$.
Let $\phi$ denote a generator of the Galois group
$\Gal(\F_{q^d}/\F_q)$. Fix a basis $\zeta_0,\dots,\zeta_{d-1}$ for
$\F_{q^d}$ over $\F_q$, where
$\zeta_i = \phi^i(\zeta_0)$.

Extend $\phi$ to an automorphism of the function field $k_1 =
\F_{q^d}(y)$ by setting $\phi(y) = y$; the fixed subfield is $k =
\F_q(y)$, of co-dimension $d$. Let $\nu_y$ denote the valuation
defined by $\nu_y(a_my^m+\dots+a_ny^n) = m$ ($a_m \neq 0$, $m<n$),
and set $F = \F_q\db{y}$, the completion with respect to $\nu_y$,
and $\O = \F_q[[y]]$, its ring of integers.

\ifXY
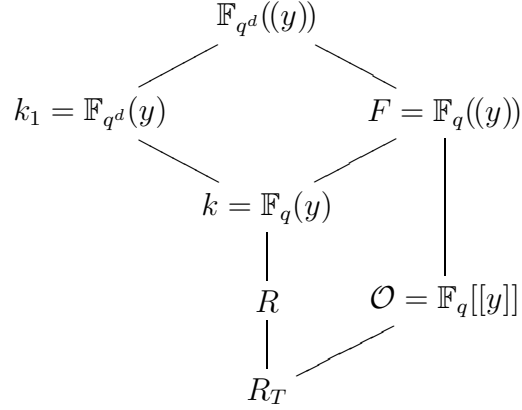
\begin{figure}
\begin{equation*}
\xymatrix@R=16pt@C=6pt{
    {}
    & {}
    & {\Fy[q^d]} \ar@{-}[dr] \ar@{-}[dl]
    & {}
    & {}
    & {}
\\
    {}
    & {k_1 = \F_{q^d}(y)} \ar@{-}[dr]
    & {}
    & {F = \Fy} \ar@{-}[dl] \ar@{-}[dd]
    & {}
    & {}
\\
    {}
    & {}
    & {k = \F_{q}(y)} \ar@{-}[d]
    & {}
    & {}
    & {}
\\
    {}
    & {}
    & {R} \ar@{-}[d]
    & {\O = \F_q[[y]]} \ar@{-}[dl]
    & {}
    & {}
\\
    {}
    & {}
    & {R_T}
    & {}
    & {}
    & {}
}
\end{equation*}
\caption{Subrings of $F_{q^d}\db{y}$}\label{fi:fields_y}
\end{figure}
\else
\bigskip
Zoo.
\bigskip
\fi

Let
$$ R = \F_q\left[y, \frac{1}{y}, \frac{1}{1+y}\right] \sub k,$$
and let $R_T$ denote the subring
\begin{equation}\label{RTdef}
R_T = \F_q\left[y,\frac{1}{1+y}\right].
\end{equation}
Since $1+y$ is invertible in $\O$, $R_T \sub \O$.

For a commutative $R_T$-algebra $S$ (namely a commutative ring
with unit which is an $R_T$-module, \eg\ $k$, $F$, $R$ or $R/I$
for an ideal $I \normali
R$), we denote by $\ww$ the element $\ww \cdot 1 \in S$.
For such $S$ we
define an $S$-algebra $\AA(S)$, by
$$\AA(S) = \bigoplus_{i,j=0}^{d-1}{S \zeta_i z^j},$$
with the relations
\begin{equation}\label{rel}
z \zeta_i = \phi(\zeta_i) z, \qquad z^d = 1+\ww.
\end{equation}

The center of $\AA(S)$ is $S$. We will frequently use the fact
that for an $R_T$-algebra $S$, $\AA(S) = \AA(R_T) \tensor[R_T] S$.
It is well known that $\AA(k)$ is a central simple
algebra, and so there is a norm map $\mul{\AA(k)} \ra \mul{k}$,
which induces a norm map $\mul{\AA(S)} \ra \mul{S}$ for every $S
\sub k$. The norm map is a homogeneous form of degree $d$ (in the
coefficients of the basis elements $\set{\zeta_i z^j}$), so the
norm is also defined for quotients $\AA(S/I)$.
We remark that $\F_{q^d}\tensor[\F_q] S$ is Galois over $S$ and
$1+\ww$ is invertible in $S$, so $\AA(S)$ is an Azumaya algebra
over $S$ (see \cite{DI}, where they are called `central separable
algebras').
This fact will not be used in the rest of the paper.

If $\AA(S) \isom \M[d](S)$, we say that $\AA(S)$ is \defin{split}.
We need a criterion for this to happen.
If $S_1 = \F_{q^d} \tensor[\F_q] S$ is a field (so necessarily $S$
is a subfield), then $\AA(S)$ is the cyclic algebra
$(S_1/S,\phi,1+\ww) = S_1[z]$ with the relations in \eq{rel}. This
is a simple algebra of degree $d$ over its center $S$. Recall
Wedderburn's norm criterion for cyclic algebras
\cite[Cor.~1.7.5]{Jac}: the algebra $\AA(S) = (S_1/S,\phi,1+\ww)$
splits iff $1+\ww$ is a norm in the field extension $S_1/S$. More
generally, the exponent of $(S_1/S,\phi,1+\ww)$ (\ie\ its order in
the Brauer group $\Br(S)$) is the minimal $i>0$ such that
$(1+\ww)^i$ is a norm. In particular, if this exponent is $d =
\dimcol{S_1}{S}$, $\AA(S)$ is a division algebra (since the
exponent of a central simple algebra is always bounded by the
degree of the underlying division algebra).

The algebra $\AA(R)$ will later be used to construct the desired
complexes. As mentioned above, $\AA(R) \sub \AA(k)$, $k$ being the
ring of fractions of $R$; moreover, $\AA(k)$ is the ring of
central fractions of this algebra,
$$\AA(k) = (R-\set{0})^{-1}\AA(R).$$

We now consider completions of $\AA(k)$. The global field $k =
\F_q(y)$ has the minus degree valuation, defined for $f,g\in
\F_q[y]$ by $\n(f/g) = \deg(g) - \deg(f)$. Also recall that the
other nonarchimedean discrete valuations of $\F_q(y)$ are in
natural correspondence with the prime polynomials of $\F_q[y]$.
For a prime polynomial $p \in \F_q[y]$, the valuation is defined
by $\nu_p(p^i f/g) = i$, where $f,g$ are polynomials prime to $p$.
The ring of $p$-adic integers in $k$ is $\F_q[y]_p = \set{f/g
\suchthat (p,g) = 1}$. The completion with respect to a valuation
$\nu = \nu_p$ is $k_{\nu} = \F_q[y]_p\db{p} =
\set{\sum_{i=-v}^{\infty}{\alpha_i p^i}}$ (which we will also
denote by $k_{p}$). The ring of integers of $k_p$ is $\O_p =
{{\F_q}[y]}_p[[p]]$.
The notation $\yy{\nu}$ is used for the degree
valuation, since the completion of $k$ with respect to this
valuation is $\F_q\db{1/y}$; moreover, filtration by the ideal
$\ideal{1/y}$ of the ring of integers $\F_q[[1/y]]$ determines the
valuation.

The Albert-Brauer-Hasse-Noether Theorem describes division
algebras over $k$ in terms of their local invariants, which
translates to an injection $\Br(k) \ra \oplus \Br(k_{p})$. More
precisely, the $d$-torsion part of  $\Br(k_{p})$ is cyclic of
order $d$ for every valuation $\nu_p$.
Taking the unramified extension $k'_{p}/k_{p}$ of dimension $d$,
these $d$ classes can be written as the cyclic
algebras $(k'_{p}/k_p,\phi,\uni^i)$ for $i = 0,\dots,d-1$, where
$\uni$ is a uniformizer.
If $(k_1/k,\phi,c)$ is a cyclic
$k$-algebra, the local invariants are determined by the values
$\nu(c)$ (\cite[Chaps.~17--18]{Pierce} is a standard reference,
though the focus is on number fields).

There are only two valuations $\nu$ of $k$ for which
$\nu(1+y) \neq 0$, namely $\nu_{1+y}$ and $\nu_{1/y}$, for which
the values are $1$ and $-1$, respectively. We thus have
\begin{prop}\label{whowhat}
The completions $\AA(\F_q\db{1/y})$ and $\AA(\F_q\db{1+y})$ of
$\AA(k)$ are division algebras. On the other hand, for any other
completion $k_\n$ of $k$, $\AA(k_\n)$ splits.
%
\end{prop}

In particular
\begin{equation}\label{ARinMF}
\AA(R) \sub \AA(k) \tensor[k] F = \AA(F) \isom \M(F).
\end{equation}
The same argument embeds
\begin{equation}\label{ARTinMO}
\AA(R_T) \sub \AA(\O) \isom \M[d](\O).
\end{equation}

\ifXY
\begin{figure}
\begin{equation*}
\xymatrix@R=24pt@C=32pt{
    {\GL[d](F)} \ar[r]^{\det}
    & {\mul{F}}
\\
    {\mul{\AA(R)}} \ar@{^(->}[u] \ar[r]^{\Norm{}}
    & {\mul{R}} \ar@{^(->}[u]
}
\end{equation*}
\caption{Norm and determinant}\label{Figdet}
\end{figure}
\else
\bigskip
Zoo 3.
\bigskip
\fi

We use the algebras $\AA(S)$ to define algebraic group schemes. For
an $R_T$-algebra $S$, let $\tG'(S) = \mul{\AA(S)}$, the invertible
elements of $\AA(S)$, and $G'(S) = \mul{\AA(S)}/\mul{S}$. Recall
that for every $R_T$-algebra $S$,
one can define a multiplicative
function $\AA(S) \ra S^\times$ called the reduced norm (\eg\ by
taking the determinant in a splitting extension of $S$). In particular,
the diagram in Figure \ref{Figdet} commutes.
We can
thus define $\tG'_1(S)$ as the set of elements of $\tG'(S)$ of
norm $1$, and $G'_1(S)$ as the image of $\tG'_1(S)$ under the map
$\tG'(S) \ra G'(S)$ (see the square in the middle of Figure
\ref{fi:GG}).

\begin{rem}\label{ker}
The sequence
$$1 \lra \mu_d(S) \lra \tG'_1(S) \lra G'_1(S) \lra 1$$
is exact, where $\mu_d(S)$ is the group of $d$-roots of unity in $S$.
\end{rem}

 These group schemes are forms of
the classical groups $\tG(S) = \GL(S)$, $G(S) = \PGL(S)$,
$\tG_1(S) = \SL(S)$ and $G_1(S) = \PSL(S)$. If $\AA(S)$ is a
matrix ring, we have that $\tG'(S) = \tG(S)$ and $G'(S) = G(S)$.

It is useful to have equivalent definitions for the groups $G'(S)$
for various rings $S$. Fix the ordered basis
$$\set{\zeta_0,\dots,\zeta_{d-1},\zeta_0z,\dots,\zeta_{d-1}z,\dots,\zeta_0 z^{d-1},\dots,\zeta_{d-1}z^{d-1}}$$
of $\AA(S)$ over $S$. Conjugation by an invertible element $a \in
\AA(S)$ is a linear transformation of the algebra. Let $i_S \co
G'(S) \ra \GL[d^2](S)$ be the induced embedding. If $S \sub S'$,
then the diagram \ifXY
\begin{equation*}
\xymatrix{
    {G'(S')} \ar@{->}[r]^{i_{S'}}
    & {\GL[d^2](S')}
\\
    {G'(S)} \ar@{^(->}[u] \ar[r]^{i_S}
    & {\GL[d^2](S)} \ar@{^(->}[u]
}
\end{equation*}
\else
\bigskip
((the inclusion $G'(S) \sub G'(S') \stackrel{i_{S'}}{\lra}
\GL[d^2](S')$ is equal to $\AA(S) \stackrel{i_{S}}{\lra}
\GL[d^2](S) \sub \GL[d^2](S')$.))
\bigskip
\fi
commutes.

\begin{prop}\label{twodefs}
Let $R_T \sub S \sub S'$ be commutative rings, such that $S$ is a
Noetherian
unique factorization domain. Then $$i_S G'(S) = i_{S'}G'(S') \cap
\GL[d^2](S),$$ the intersection taken in $\GL[d^2](S')$.
\end{prop}
\begin{proof}
The inclusion $i_S G'(S) \sub i_{S'}G'(S') \cap \GL[d^2](S)$ is
trivial since $\AA(S) \sub \AA(S')$. Let $\alpha \in i_{S'}G'(S')
\cap \GL[d^2](S)$, then $\alpha$ is an isomorphism of algebras
(since it is induced by an element of $\AA(S')$) and preserves $S$
(as it belongs to $\GL[d^2](S)$). It is thus an automorphism of
$\AA(S)$, which must be inner \cite[Thm.~3.6]{AG}.
\end{proof}

The proposition covers, in particular, $S = R_T, R, k, k_{\nu},
\O_{\nu}$, as well as $\bar{R}_T$, and $\bar{R}$ which are defined
in Section \ref{splitx}, with an arbitrary extension $S'$ (usually
taken to be from the same list).

\begin{prop}
$G'(R)$ is a discrete subgroup of $G(F)$.
\end{prop}
\begin{proof}
The ring $R = \F_q[\frac{1}{\ww},\ww,\frac{1}{1+\ww}]$ embeds
(diagonally) as a discrete subgroup of the product
$$
F \times \yy{k} \times k_{p}.
$$
This can be seen by letting $a_n =
\frac{f_n(\ww)}{\ww^{i_n}(1+\ww)^{j_n}}$ be a sequence of non-zero
elements in $R$ (with $f_n(\lam) \in \F_q[\lam]$ and $i_n, j_n
\geq 0$) such that $a_n \ra 0$.  We then have that $\nu_y(a_n) \ra
\infty$, which implies $i_n = 0$ for $n$ large enough and
$\nu_y(f_n(\ww)) \ra \infty$. Likewise, $\nu_p(a_n) \ra \infty$,
so $j_n = 0$ for $n$ large enough.  This implies that
$\yy{\nu}(a_n) \ra -\infty$, which is a contradiction.

It follows that the diagonal embedding of $G'(R)$ into
$$
G'(F) \times G'(\yy{k}) \times G'(k_{p})
$$
is discrete. But an algebraic group over a local field is compact
iff it has rank zero \cite{PR}.  Therefore, by \Pref{whowhat},
$G'(\yy{k})$ and  $G'(k_{p})$ are compact, and $G'(R)$ is discrete
in the other component $G'(F) = G(F)$.
\end{proof}

In fact, from general results it follows
that $G'(R)$ is a cocompact lattice
in $G'(F)$, but we will show it directly when we demonstrate that
$G'(R)$ acts transitively on the vertices of the affine building
of $G(F) = \PGL(F)$.

Consider $G(\O) = \PGL[d](\O)$, a maximal compact subgroup of
$G(F)$, which is equal to $G'(\O)$ by \Eq{ARTinMO}. Viewing $K =
i_{\O} G'(\O)$ and $i_R G'(R)$ as subgroups of $i_F G'(F) \sub
\GL[d^2](F)$, the intersection
\begin{equation}\label{zinK}
i_{R_T} G'(R_T) = K \cap i_R G'(R)
\end{equation}
is finite, being the intersection of discrete and compact
subgroups (note that $R \cap \O = R_T$).

\begin{prop}\label{GRT=1}
$G'(R_T)$ is a semidirect product of $\sg{z} \isom \Z/d\Z$ acting
on $\mul{\F_{q^d}}/\mul{\F_{q}}$.
\end{prop}
\begin{proof}
Recall that $R_T = \F_q[\ww,1/(1+\ww)]$, so that $\AA(R_T) =
\F_{q^d}[\ww,1/(1+\ww),z]$ with the relations $z\alpha z^{-1} =
\phi(\alpha)$ ($\alpha \in \F_{q^d}$) and $z^d =1+\ww$. Setting
$\ww = z^d-1$, we see that $\AA(R_T) = \F_{q^d}[z,z^{-1}]$ is a
skew polynomial ring with one invertible variable over $\F_{q^d}$.
Every element of $\AA(R_T)$ has a monomial $\alpha z^r$ ($\alpha \in
\mul{\F_{q^d}}$) with $r$ maximal, called the upper monomial (with
respect to $z$), and similarly every element has a lower monomial.
The upper monomial of a product $fg$ is equal to the product of
the respective upper
monomial, and likewise for the lower monomials.

Now let $f,g \in \AA(R_T)$ be elements with $fg = 1$, then the
product of the upper monomials and that of the lower monomials are both
equal to $1$,
proving that $f$ and $g$ are monomials. Thus, the invertible
elements of $\AA(R_T)$ are $\tG(R_T) = \set{\alpha z^{i} \suchthat
\alpha \in \mul{\F_{q^d}},\, i\in \Z}$. The result is obtained by
taking this modulo the center.
\end{proof}

\section{A simply transitive action on $\B$}\label{sec:gamma}

We continue with the notation of the last section. The embedding
\eq{ARinMF} of $\AA(R)$ into $\M(F)$ induces embeddings $\tG'(R)
\hra \GL(F)$ and $G'(R) \hra \PGL(F)$, so $G'(R)$ acts on $\B^0$.
Notice that $G'(k)$ is dense in $G'(F) = G(F)$, so its action on
$\B^0$ is transitive.

In \cite{CS1}, Cartwright and Steger present a subgroup of
$G'(k)$, which acts simply transitively on the vertices of the
building $\B$. Identifying $\B^0$ with a group, when possible, is
an important tool in the description of finite quotients of $\B$.
We will use the Cartwright-Steger group, constructed in a
different way. In particular we construct the group as a subgroup
of $G'(R)$, which is easily shown to be discrete.
The proof that it acts transitively relies on the existence of an element
of the appropriate norm, and that the action is simple is seen by
a relatively easy computation of the group scheme over various rings.

Recall that $\AA(R) = R[\zeta_i,z]$ where $z\zeta_i z^{-1} =
\phi(\zeta_i)$ and $z^d = 1+\ww$, and let
$$b = 1-z^{-1} \in \AA(R).
$$

\begin{prop}\label{normb}
The reduced norm of $b$ in $\AA(R)$ is $\ww/(1+\ww)$.
\end{prop}
\begin{proof}
Recall that $\AA(R) \sub \AA(k)$ and that the restriction of the
norm function of $\AA(k)$ to $\AA(R)$ is the norm function of
$\AA(R)$. As is the case for matrices, if $a \in \AA(R)$ generates
a subalgebra of dimension $d$, its norm $\Norm(a)$ is the free
coefficient in the characteristic polynomial of $a$ multiplied by
$(-1)^d$.

The minimal polynomial of $z$ is $\lam^d - (1+\ww)$, so that
$\Norm(z) = (-1)^{d-1}(1+\ww)$. Likewise, the minimal polynomial
of $z-1$ is $(\lam+1)^d-(1+\ww)$, with the free coefficient $-y$,
so that $\Norm(z-1) = (-1)^{d-1}y$. Finally, since $b =
z^{-1}(z-1)$, $\Norm(b) = \Norm(z-1)/\Norm(z) = \ww/(1+\ww)$.
\end{proof}

Since $1+\ww$ is invertible in $\O$, we claim
\begin{cor}\label{Nb=y}
$\Norm(b) \equiv y \pmod{\mul{\O}}$.
\end{cor}

Under the embedding \eq{ARTinMO}, $b \in \AA(R_T) \sub \M(\O)$, so
$b L_0 \sub L_0$, where $L_0 = \O^d$ is the standard lattice. By
the corollary, and since the diagram in Figure \ref{Figdet}
commutes, we have $\dimcol{L_0}{b L_0} = q$. Note that $b$ is not
invertible in $\AA(R_T)$ (since the norm $\ww/(1+\ww)$ is not
invertible in $R_T$), but $b$ is invertible in $\AA(R)$.

Let $\Omega$ denote the set of neighbors of $[L_0]$ in $\B$ which
have color $1$; they correspond to sublattices of index $q$ in
$L_0$. Consider the element
$$b_u = ubu^{-1} =
1-\frac{u}{\phi(u)}z^{-1}$$
where $u \in \mul{\F_{q^d}} \sub \mul{\AA(R)}$. The embedding
$\AA(R) \sub \M(F)$ extends the regular embedding $\F_{q^d} \hra
M_d(\F_q)$ (via the basis $\zeta_0,\dots,\zeta_{d-1}$), and in
particular it takes $u$ to a matrix with coefficients in $\F_q$.
This proves that $uL_0 = L_0$, and since $\GL(F)$ preserves the
structure of $\B$, $u$ permutes the vertices in $\Omega$.
Moreover, the sublattices of index $q$ of $L_0$ correspond to the
points of the projective space $\mul{\F_{q^d}}/\mul{\F_{q}}$ via
the isomorphism $L_0/y L_0 \isom \F_{q^d}$. The action of $u$ on
the latter space is by multiplication, so $\mul{\F_{q^d}} \sub
\mul{\AA(R)}$ acts transitively on $\Omega$. Finally, we have that
$b_u (L_0) = u b u^{-1} L_0 = u(b L_0)$, so we proved:
\begin{prop}
For every lattice $L_1 \in \Omega$, there is an element $u \in
\mul{\F_{q^d}}$ (unique up to multiplication by $\mul{\F_q}$) such
that $b_u L_0 = L_1$.
\end{prop}

\begin{defn}\label{Gammadef}
Let $\tilde{\Gamma}$ denote the subgroup
\begin{equation}\label{Gamtil}
\tilde{\Gamma} = \sg{b_u \suchthat u \in \mul{\F_{q^d}}/\mul{\F_{q}}}
\end{equation}
of $\tG'(R) = \mul{\AA(R)}$, and let $\Gamma$ be the image of
$\tilde{\Gamma}$ under the projection $\tG'(R) \ra G'(R)$.
\end{defn}

\begin{prop}\label{trans}
$\Gamma$ acts transitively on $\B^0$.
\end{prop}
\begin{proof}
We show that $L_0$ can be taken to any of its sublattices $L$ by
an element of $\Gamma$. The proof is by induction on
$\dimcol{L_0}{L}$. If $\dimcol{L_0}{L} = q$, then $L \in \Omega$
and $b_u L_0 = L$ for some $u$. Assume $\dimcol{L_0}{L} > q$, and
let $L \sub L' \sub L_0$ be an intermediate sublattice such that
$\dimcol{L'}{L} = q$. By the induction hypothesis there is an
element $c \in \Gamma$ such that $c L_0 = L'$. The lattices in
$c\Omega$ are the sublattices of index $q$ of $L'$, so for some $u
\in \mul{\F_{q^d}}$, we have that $c (b_u L_0) = L$, which proves
our claim.
\end{proof}

\begin{cor}\label{actk}
If $L \sub L_0$ and $\dimcol{L_0}{L} = q^\ell$, then there are
$u_1,\dots,u_{\ell} \in \mul{\F_{q^d}}/\mul{\F_{q}}$ such that
$[L] = b_{u_1} \ldots b_{u_{\ell}} [L_0]$.

Moreover, if $L \sub L' \sub L_0$ where $\dimcol{L_0}{L'} = q^k$
and $[L'] = b_{u_1} \ldots b_{u_k} [L_0]$, then for suitable
$u_{k+1},\dots,u_{\ell}$ we have $[L] = b_{u_1} \ldots b_{u_k}
b_{u_{k+1}} \ldots b_{u_{\ell}} [L_0]$.
\end{cor}

Let $R_0 = \F_q[1/\ww]$. As $R_0$ does not contain $R_T$,
$\AA(R_0)$ is not defined, so we cannot define $G'(R_0)$ in the
usual way. However, following \Pref{twodefs}, we set $G'(R_0) =
i_{k}G'(k) \cap \GL[d^2](R_0)$.

As mentioned in the introduction, Cartwright and Steger present in
\cite{CS1} a group $\Gamma'$ which acts simply transitively on
$\B^0$. Their group is defined as follows:
\begin{defn}\label{Gammaprimedef}
The group $\Gamma'$ is composed of the elements of $G'(R_0)$ which
modulo $1/\ww$ are upper triangular with identity $d\times d$
blocks on the diagonal.
\end{defn}

A careful computation, using $b_u^{-1} =
1+\frac{1}{\ww}\sum_{i=0}^{d-1}{ u\phi^{i}(u)^{-1}z^{i}}$, reveals
that
\begin{equation}\label{buis}
b_u (\zeta z^k) b_u^{-1} = \zeta z^k + f_{u,k}(\zeta)
\sum_{i=0}^{k-1}\frac{u}{\phi^{i}(u)}z^{i} + \frac{1}{\ww}
f_{u,k}(\zeta) \sum_{i=0}^{d-1}\frac{u}{\phi^{i}(u)}z^{i}
\end{equation}
for every $0 \leq k < d$ and $\zeta \in \F_{q^d}$, where
$f_{u,k}(\zeta) =
\frac{\phi^{k}u}{u}\zeta-\frac{\phi^{k-1}u}{\phi^{-1}u}\phi^{-1}\zeta$.

It follows that $b_u \in \Gamma'$, so $\Gamma = \sg{b_u} \sub
\Gamma'$. This leads to an easy proof of the main property of
$\Gamma$:
\begin{prop}\label{simply}
The group $\Gamma$ acts simply transitively on $\B^0$, namely, for
$K=G(\O)$,
\begin{eqnarray}
\Gamma \cdot K &= &G(F), \label{GammaK=G} \\
\Gamma \cap K &= &1. \label{GammaK=1}
\end{eqnarray}
In particular, $\Gamma = \Gamma'$.
\end{prop}
\begin{proof}
The first equation, equivalent to $\Gamma$ acting transitively on
$\B^0 = G(F)/K$, was proved in \Pref{trans}. To compute the
intersection, note that $\Gamma \cap K \sub \Gamma' \cap K \sub
G'(R) \cap K = G'(R_T)$ by \Eq{zinK}, so $\Gamma \cap K \sub
G'(R_T) \cap \Gamma'$. Let $\zeta z^j$ be an element of $G'(R_T)
\cap \Gamma'$ for $\zeta \in \mul{\F_{q^d}}$ and $j = 0,\dots,d-1$
(see \Pref{GRT=1}). By the definition of $\Gamma'$, the $d^2\times
d^2$ matrix representing $\zeta z^j$ should be upper triangular
modulo $1/\ww$. Now, $(\zeta z^j) z (\zeta z^j)^{-1} =
\frac{\zeta}{\phi(\zeta)} z$, so $\zeta \in \mul{\F_q}$.
Furthermore for every $a \in \F_{q^d}$, $(\zeta z^j) a (\zeta
z^j)^{-1} = \phi^j(a)$, proving that $j = 0$ and $\zeta z^j$ is
central; thus $\Gamma' \cap K = 1$. This proves that $\Gamma'$
acts simply transitively on $G(F)/K$, and since $\Gamma \sub
\Gamma'$, these groups are equal.
\end{proof}

\begin{prop}\label{GRGR0}
The group $\Gamma$ is a normal subgroup of $G'(R)$, and $G'(R) =
G'(R_T) \semidirect \Gamma$. Moreover,
\begin{equation}\label{11}
G'(R_0) = G'(R).
\end{equation}
\end{prop}
\begin{proof}
Let $g \in G'(R)$, and write $g = \gamma a$ for $\gamma \in \Gamma
\sub G'(R)$ and $a \in K$. Then $a = \gamma^{-1} g \in K \cap
G'(R) = G'(R_T)$, by \Eq{zinK}, and $G'(R) = \Gamma \cdot
G'(R_T)$.

The generators of $\Gamma$ are permuted by $G'(R_T)$: $(\zeta_i
z^j) b_u (\zeta_i z^j)^{-1} = b_{\zeta_i \phi^j(u)}$, so $\Gamma
\normali G'(R)$. Moreover, $\Gamma \cap G'(R_T) \sub \Gamma \cap K
= 1$, proving the second claim.

Finally, $G'(R_0) \sub G'(R)$ but direct inspection using
\Eq{buis}, shows that $\Gamma$ and $G'(R_T)$ are contained in
$\GL[d^2](R_0)$ (see \cite[Thm.~2.6]{CS1}).
\end{proof}

Since $\Gamma$ acts simply transitively, we can identify $\B^0$
with $\Gamma$ by $\gamma \mapsto \gamma[L_0]$. Recall that by
\Eq{ARinMF}, $\Gamma \sub \PGL(F)$. The determinant (modulo
scalars) thus takes $\gamma \mapsto \det(\gamma)\cdot {\mul{F}}^d
\in \md{F}$. Considering this modulo $\mul{\O}$ (which is
superfluous if $d$ is prime to $q$), we obtain the cyclic group
$\md{F}\mul{\O} \isom \sg{y}/\sg{y^d}$. Moreover, since
$\color(b_u \chix) = \color(\chix) + 1$ and $\Norm(b_u) \equiv y
\pmod{\mul{\O}}$ by \Cref{Nb=y}, the diagram in Figure
\ref{detcol} commutes. We let $\delta$ denote the color map
$\Gamma \ra \Z/d$. Letting $\Gamma_1 = \Ker(\delta)$, the short
sequence
$$1 \lra \Gamma_1 \lra \Gamma \stackrel{\delta}{\lra} \Z/d\Z \lra 1$$
is exact.
\begin{cor}\label{Gamma1} $\dimcol{\Gamma}{\Gamma_1}
= d$.
\end{cor}

\ifXY
\begin{figure}
\begin{equation*}
\xymatrix@R=24pt@C=32pt{
    {\Gamma} \ar@{->}[rr]^(0.35){\det} \ar[drr]^{\delta}
    \ar[d]_{\gamma\mapsto \gamma[L_0]}
    &
    & {\md{F}}\mul{\O}  \ar[d]^{y \mapsto 1}
\\
    {B_0} \ar[rr]^{\color}
    &
    & {\Z/d\Z}
}
\end{equation*}
\caption{Determinant and coloring}\label{detcol}
\end{figure}
\else
\bigskip
Zoo 3.
\bigskip
\fi

\section{The geometry of $\Gamma$}\label{geom}

By definition, $\Gamma$ is generated by the elements $b_u =
u(1-z^{-1})u^{-1}$ of $G'(R)$, where $u$ ranges over
$\mul{\F_{q^d}}/\mul{\F_{q}}$.

Since $b_u = 1 - \frac{u}{\phi(u)}z^{-1}$, for every $u$, there is
a unique $r \in \mul{\F_{q^d}}$ with $\Norm(r)=1$, such that $b_u
= b_{(r)} = 1 - rz^{-1}$, where $\Norm =\Norm_{\F_{q^d}/\F_q}$ is
the norm map.
\begin{rem}\label{uvprime}
For every $u\neq v \in \mul{\F_{q^d}}/\mul{\F_{q}}$ there are unique
$u',v' \in \mul{\F_{q^d}}/\mul{\F_{q}}$ such that
\begin{equation}\label{pairs}
b_{u} b_{u'} = b_{v} b_{v'}.
\end{equation}
\end{rem}
\begin{proof}
The conceptual reason is that since $b_u L_0 \neq b_v L_0$, their
intersection $L = b_u L_0 \cap b_v L_0$ is a submodule of index
$q^2$ of $L_0$ (which contains $yL_0$), and so by \Cref{actk}
there are unique $u',v'$ such that $L = b_{u}b_{u'} L_0 =
b_{v}b_{v'}L_0$. Since $\Gamma$ acts simply transitively,
$b_{u}b_{u'} = b_{v}b_{v'}$.

Computationally, assume first $d > 2$. Since
$(1-rz^{-1})(1-r'z^{-1}) = 1 - (r+r')z^{-1} +
r\phi^{-1}(r')z^{-2}$, $b_{(r)}b_{(r')} = b_{(s)}b_{(s')}$ in
$\Gamma$ iff
$$r+r' = s+s'\quad\,\mbox{and}\quad\, r\phi^{-1}(r') = s\phi^{-1}(s').$$
The unique solution is $r' = \frac{s-r}{\phi s -\phi r} \phi(s)$
and $s' = \frac{s-r}{\phi s -\phi r} \phi(r)$, where $r'$ and $s'$
indeed have norm $1$.

If $d = 2$ then (as elements of $\tilde{\Gamma}$) $b_{(-r)}b_{(r)}
= 1 - \frac{r\phi(r)}{1+y}$, which is in the center $R =
\Cent(\AA(R))$. Therefore in $\Gamma$, $b_{(r)}^{-1} = b_{(-r)}$.
Now the unique solution to the above equation is $b_{(r)}b_{(-r)}
= 1 = b_{(s)}b_{(-s)}$.
\end{proof}

The following is essentially a reformulation of the simply
connectedness of $\B$.
\begin{thm}\label{presGamma}
Let $\Gamma'$ be the abstract group generated by $\set{x_u
\suchthat u \in \mul{\F_{q^d}}/\mul{\F_{q}}}$, with the relations
$x_{u'} x_{u} = x_{v'} x_{v}$ whenever $b_{u'} b_{u} = b_{v'}
b_{v}$ in $\Gamma$, and a single relation $x_{u_1} \ldots x_{u_d}
= 1$ for some $u_1,\dots,u_d$ such that
\begin{equation}\label{yword}
r = b_{u_1} \ldots b_{u_d} = 1
\end{equation}
in $\Gamma$.
Then $x_u \mapsto b_u$ is an isomorphism $\Gamma' \ra \Gamma$.
\end{thm}
\begin{proof}
Since $\Gamma$ acts simply transitively, $b_{u_1}^{\epsilon_1}
\ldots b_{u_n}^{\epsilon_n} = 1$ iff $b_{u_1}^{\epsilon_1} \ldots
b_{u_n}^{\epsilon_n} L = y^i L$ for some lattice $L$ and $i \in
\Z$.

\begin{figure}[!h]
\ifXY
\begin{equation}\nonumber
\xymatrix@R=12pt{
    {}
    & L_0  \ar@{-}[rd] \ar@{-}[ld]
    & {}
\\
    b_{u_1}L_0 \ar@{-}[rd] \ar@{-}@/_20pt/[rdddd]
    & {}
    & b_{v_1}L_0 \ar@{-}[ld] \ar@{-}@/^20pt/[ldddd]
\\
    {}
    & \delta L_0 \ar@{-}[ddd]
    & {}
\\
    {}
    & {}
    & {}
\\
    {}
    & {}
    & {}
\\
   {}
    & \gamma L_0
    &
}
\end{equation}
\else

    {}
    & [L_0]
    & {}
\\
    {}
    & {}
    & {}
\\
    {}
    & {}
    & {}
\\
    {}
    & \circ \ar@{-}[uuu]
    & {}
\\
    \circ \ar@{-}[ru]_{v'} \ar@{-}@/^20pt/[ruuuu]^{u_2\ldots u_n}
    & {}
    & \circ \ar@{-}[lu]^{u'} \ar@{-}@/_20pt/[luuuu]_{v_2\ldots v_n}
\\
   {}
    & [L_0] \ar@{-}[ru]_{v_1} \ar@{-}[lu]^{u_1}
    &

\bigskip

(Diagram: going down in three ways).

\bigskip
\fi
\caption{Proof of \Tref{presGamma}}\label{pic1}
\end{figure}

We first claim that every equation $\gamma = b_{u_1} \ldots
b_{u_n} = b_{v_1} \ldots b_{v_n}$ (both expressions of the same
length) can be obtained from the relations of type \eq{pairs}.
Indeed, for $n \leq 2$ the claim is obvious.

If $u_n = v_n$ then we are done by induction, so we assume $u_n
\neq v_{n}$. By the above remark there are $u'$ and $v'$ such that
$\delta = b_{u_n}b_{u'} = b_{v_n}b_{v'}$, and then $\delta L_0 =
b_{u_1}L_0 \cap b_{v_1}L_0$ is of index $q^2$ in $L_0$, and
contains $\gamma L_0$; see Figure \ref{pic1}. By \Cref{actk} there
are $w_3,\dots,w_n$ such that $\gamma L_0 = \delta b_{w_3} \ldots
b_{w_n} L_0$. By induction the equations $b_{u_2} \ldots b_{u_n} =
b_{u'} b_{w_3} \ldots b_{w_n}$ and $b_{v_2} \ldots b_{v_n} =
b_{v'} b_{w_3} \ldots b_{w_n}$ follow from the relations
\eq{pairs}, and so is the case with $b_{u_1} b_{u_2} \ldots
b_{u_n} = b_{u_1} b_{u'} b_{w_3} \ldots b_{w_n} = b_{v_1} b_{v'}
b_{w_3} \ldots b_{w_n} = b_{v_1} b_{v_2} \ldots b_{v_n}$, which
proves our claim. In particular, all the relations of the form
\eq{yword} are equivalent modulo the relations \eq{pairs}.

Let $u \in \mul{\F_{q^d}}/\mul{\F_{q}}$. Using \Cref{actk} we  see
that $b_{u}^{-1}$ is equal to a product of $d-1$ generators. Given
a relation, we can thus assume it has the form $b_{u_1} \ldots
b_{u_n} = 1$. Then $b_{u_1} \ldots b_{u_n} L_0 = y^i L_0$ where $n
= di$, so $b_{u_1} \ldots b_{u_n} = r^i$ is an equality of
elements of equal length, which follows from the relations
\eq{pairs} by the above claim.

\end{proof}

\begin{rem}\label{reld}
If $d$ is a power of $\mychar F$, then equation \eq{yword} takes a
particularly simple form: $b_1^d = (1-z^{-1})^d = 1-z^{-d} =
\frac{\ww}{1+\ww} \in \mul{R}$, so $b_1^d = 1$ in $\Gamma$.
\end{rem}

{}From the theorem, together with the case $d = 2$ of Remark
\ref{uvprime}, we obtain
\begin{cor} If $d = 2$, then the defining relations of $\Gamma$
are $b_{(r)}^{-1} = b_{(-r)}$. In particular, if $q$ is odd then
$\Gamma$ is free on $\frac{q+1}{2}$ generators, and if $q$ is even
then $\Gamma$ is a free product of $q+1$ cyclic groups of order
$2$.
\end{cor}

\section{Finite quotients of $\B$}\label{sec:finq}

In this section we specify quotients of $\Gamma$, whose Cayley
graphs (with respect to the generators $b_u$) define Ramanujan
complexes. These quotients are shown to be subgroups of $\PGL$
which contain $\PSL$, over a pre-determined finite local ring.

Recall the definitions of $\tG', \tG'_1, G'$ and $G'_1$ from
Section \ref{sec2}. Let $S$ be an $R_T$-algebra, and $I \normali
S$ an ideal.  The map $S \ra S/I$ induces an epimorphism of
algebras $\AA(S) \ra \AA(S/I)$. It also induces a homomorphism
$\tG'(S) \ra \tG'(S/I)$ (which, in general, is not onto). The
kernel of this map is called the congruence subgroup of $\tG'(S)$
with respect to $I$ and denoted by $\con{\tG'}{S}{I}$. In a
similar manner we have a map $G'(S) \ra G'(S/I)$ defined by $a
\mul{S} \mapsto (a+I \AA(S))\mul{(S/I)}$, and the kernel is
denoted $\con{G'}{S}{I}$. Note that $\con{\tG}{S}{I}$ is mapped
into $\con{G}{S}{I}$, but this map is not necessarily onto. We set
$\con{\tG'_1}{S}{I} = \tG'_1(S) \cap \con{\tG'}{S}{I}$, and
likewise $\con{G'_1}{S}{I} = \con{G'}{S}{I} \cap G'_1(S)$. (Again,
note that $\con{G'_1}{S}{I} = \set{a\mul{S}\suchthat
a\in\tG'_1(S),\,a\in \mul{S}+I\AA(S)}$, while $\con{\tG'_1}{S}{I}
= \set{a \suchthat a\in\tG'_1(S),\,a \in 1+I\AA(S)}$, so the map
$$\con{\tG'_1}{S}{I} \ra \con{G'_1}{S}{I}$$ induced by $\tG'(S)
\ra G'(S)$ may not be onto). See Figure \ref{fi:GG}.

\ifXY
\begin{figure}[!h]
\begin{equation*}
\xymatrix@C=20pt{
    {}
    & {\tG'(S,I)} \ar@{^(->}[rr] \ar[dd]|[d]\hole
    & {}
    & {\tG'(S)} \ar@{->}[rr] \ar[dd]|[d]\hole
    & {}
    & {\tG'(S/I)} \ar@{->}[dd]
\\
    {\con{\tG'_1}{S}{I}} \ar@{^(->}[rr] \ar@{^(->}[ru] \ar[dd]
    & {}
    & {\tG'_1(S)} \ar[dd] \ar@{^(->}[ru] \ar[rr]
    & {}
    & {\tG'_1(S/I)} \ar[dd] \ar@{^(->}[ru]
    & {}
\\
    {}
    & {G'(S,I)} \ar@{^(->}[rr]|(0.525)\hole
    & {}
    & {G'(S)} \ar[rr]|(0.477)\hole
    & {}
    & {G'(S/I)}
\\
    {\con{G'_1}{S}{I}} \ar@{^(->}[rr] \ar@{^(->}[ru]
    & {}
    & {G'_1(S)} \ar@{^(->}[ru] \ar[rr]
    & {}
    & {G'_1(S/I)} \ar@{^(->}[ru]
    & {}
}
\end{equation*}
\caption{Group relations of Section \ref{sec:finq}}\label{fi:GG}
\end{figure}
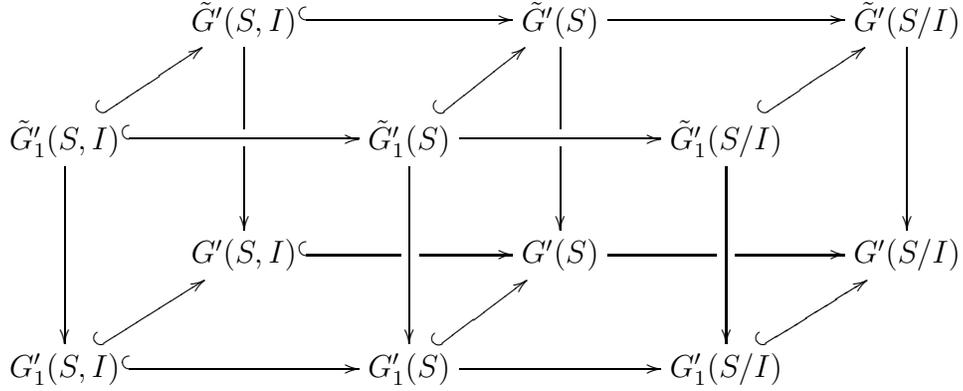
\else
\bigskip
Zoo 4.
\bigskip
\fi

\medskip

Let $I \normali R$. We define the congruence subgroup of $\Gamma$ to be
\begin{equation}\label{defGI}
\Gamma(I) = \Gamma \cap \con{G'}{R}{I},
\end{equation}
a normal subgroup of $G'(R)$ (by \Pref{GRGR0}). Recall that by
\Pref{simply}, $\Gamma$ acts simply transitively on the vertices
of $\B = G(F)/K$, where $K = G(\O)$. The set of vertices in the
quotient complex $\B_I = \dc{\Gamma(I)}{G(F)}{K}$ is isomorphic to
$\Gamma/\Gamma(I)$. To make this an isomorphism of complexes, we
need to define a complex structure on $\Gamma/\Gamma(I)$.

The generators $\set{b_u}$ correspond to the neighbors of color
$1$ of $[L_0]$. More generally, by \Cref{actk}
the neighbors of color $k$ correspond to the products $b_{u_1}
\ldots b_{u_k}$ which are headers of a product of $d$ generators
which equals $1$. The $1$-skeleton of $\B_I$ is thus the Cayley
graph of $\Gamma/\Gamma(I)$ with respect to these generators (for
$k = 1,\dots,d-1$), and there are `partial Laplacian operators'
$A_k$ of $\Gamma/\Gamma(I)$, induced from the Hecke operators of
$\B$. Frequently, these are the colored Laplacian operators (\ie\
they can be defined for $\Gamma/\Gamma(I)$ directly); we address
this issue towards the end of this section. The higher dimensional
cells of $\Gamma/\Gamma(I)$ are then defined to make it a clique
complex. Thus, $\dc{\Gamma(I)}{G(F)}{K}$ and $\Gamma/\Gamma(I)$
become isomorphic complexes.

\begin{rem}\label{NGam}
If $f(\lam) \in \F_q[\lam]$ is divisible by $\lam$, then
$\con{G'}{R_0}{\ideal{f(1/y)}} \sub \con{G'}{R_0}{\ideal{1/y}}
\sub \Gamma$, by \Dref{Gammaprimedef} of $\Gamma'$ and \Eq{11}. In
this case $\Gamma(I) = \con{G'}{R_0}{I}$ (where $I =
\ideal{f(1/y)}$) so $\dom{\Gamma(I)}{\B} \isom \Gamma /
\Gamma(I)$.
\end{rem}

 Let $C = \set{(z_1,\dots,z_d)\in \C^d \suchthat \abs[]{z_i}
=1,\, z_1 \ldots z_d = 1}$ and $\sigma \co C \ra \C^{d-1}$ the map
defined by $(z_1,\dots,z_d) \mapsto (\lam_1,\dots,\lam_{d-1})$,
where
$$\lam_k = q^{k(d-k)/2}\s_k(z_1,\dots,z_d)$$
and $\s_k$ are the symmetric functions. Then $\Aspec_d =
\sigma(C)$ is the simultaneous spectrum of the Hecke operators
$A_1,\dots,A_{d-1}$ acting on $\B$ \cite[Thm.~2.9]{paperI}.

Recall \cite{paperI} that $\B_I$ is called a Ramanujan complex if
the eigenvalues of every non-trivial simultaneous eigenvector $v
\in \L2{\B_I}$, $A_k v = \lam_k v$, satisfy
$(\lam_1,\dots,\lam_{d-1}) \in \Aspec_d$. It can be shown that the
closure of the union of the sets of simultaneous eigenvalues over
any family of quotients of $\B$ (with unbounded injective radius)
contains $\Aspec_d$.

The following theorem relies on one of the main results of
\cite{paperI}, for which we assume the global Jacquet-Langlands
correspondence for function fields. See the introduction for more
details.

\begin{thm}\label{main}
For every $d\geq 2$ and every $0 \neq I \normali R$,
the Cayley complex of $\Gamma / \Gamma(I)$
is a Ramanujan complex.
\end{thm}
\begin{proof}
Write $I = R \pp$ for some polynomial $\pp \in \F_q[y]$ which is
prime to $y$ and $1+y$, then $I \cap R_0$ is generated by
$y^{-e}\pp$ where $e = \deg(\pp)$. Let $I' = \ideal{y^{-e-1}\pp}
\normali R_0$. Set $N = \con{G'}{R_0}{I'}$, then $N \sub \Gamma$
by \Rref{NGam}, and hence $N \sub \Gamma(I)$ (this property is not
satisfied by the more `natural' candidate
$\con{G'}{R_0}{\ideal{y^{-e}\pp}}$).

In the notation of \cite{paperI}, $D = \AA(k)$ is the division algebra used to
define $G' = \mul{D}/\mul{Z}$, $F$ is the completion of $k$ with respect to
$\nu_0 = \nu_y$, and the ramification places
are $T = \set{\yy{\nu},\nu_{1+y}}$. The ring $R_0$ is
\begin{equation*}
R_0 = \bigcap_{\nu \in \Vals \minusset \set{\valF}} (k \cap
\O_\nu) = \F_q[1/y],
\end{equation*}
where $\Vals$ is the set of  valuations of $k$. Since
$$R = \bigcap_{\nu \in \Vals \minusset T \cup \set{\valF}} (k \cap \O_\nu),$$
every ideal of $R$ restricts to an ideal of $R_0$ which is prime to $T$.

Note that $I'$ is prime to $\nu_{1+y}$ (since $\nu_{1+y}(\pp) =
0$). In \cite[Thm.~1.3(a)]{paperI} we prove that $X = \dom{N}{\B}
= \dc{N}{\Gamma K}{K}$ is Ramanujan. This complex is isomorphic to
$\Gamma / N$ by $N \gamma K \mapsto \gamma N$, so $\Gamma / N$ is
Ramanujan.
Since $N \sub \Gamma(I)$, $\Gamma / \Gamma(I)$ is also Ramanujan.
\end{proof}

Our next goal is to identify $\Gamma/\Gamma(I)$ as an abstract
group. We assume that $I$ is a power of a prime ideal (so $R/I$ is
a finite local ring). Let $\ideal{\pp} = \sqrt{I}$ denote the
radical of $I$, so $I = \ideal{\pp^s}$ for some $s \geq 1$, where
$\pp \in \F_q[y]$ is prime. By definition of $\Gamma(I)$, we have
that
\begin{eqnarray*}
\Gamma/\Gamma(I) & = & \Gamma/(\Gamma \cap \con{G'}{R}{I}) \\
    & \isom & (\Gamma\cdot \con{G'}{R}{I})/\con{G'}{R}{I} \\
    & \sub & G'(R)/\con{G'}{R}{I} \\
    & \hookrightarrow & G'(R/I) \\
    & \isom & \PGL(R/I),
\end{eqnarray*}
where the final isomorphism follows from the definition of
$G'(R/I)$ as $\mul{\AA(R/I)}$ modulo center, and Wedderburn's
theorem which implies that $\AA(R/I) \isom \M[d](R/I)$.

\begin{thm}\label{isonto}
The map $G'_1(R) \ra G'_1(R/I)$ is onto.
\end{thm}
\def\GGG{\tG'_1}
\begin{proof}
Let $\A$ denote the ring of ad\`{e}les over $k$, and $\A_0 = \A/F$.
Since $\GGG$ is connected, simply connected and simple, and
$\GGG(F)$ is non-compact, the strong approximation theorem
(\cite[Thm.~7.2]{PR}, \cite{Prasad}) asserts that $\GGG(k)\GGG(F)$
is dense in $\GGG(\A)$ (see \cite[Sec.~3.2]{paperI}).
Therefore, $\GGG(k)$ is dense in $\GGG(\A_0)$. Let $U =
\GGG(k_{1/y})\GGG(k_{1+y})\prod_{\nu \neq
\yy{\nu},\nu_{1+y},\nu_y} \GGG(\O_{\nu})$, an open subgroup of
$\GGG(\A_0)$. Then $\GGG(R) = \GGG(k) \cap U$ is dense in $U$.

In particular, $G_1'(R)$ is dense in $G_1'(\O_\pp)$, where
$\O_\pp$ is the ring of integers in the completion $k_{\pp}$.
 (Notice that $R/I \isom \O_\pp/P_\pp^s$ where $P_\pp$ is the
maximal ideal of $\O_\pp$). Since $G_1'(R/I) \isom \PSL(R/I)$
which is finite, it is enough to show that the map from
$G_1'(\O_\pp)$ to $G_1'(R/I)  = G_1'(\O_\pp/P_\pp^s)$ is onto.

Let $c \in G_1'(\O_\pp/P_\pp^s)$, so $c$ is an automorphism of
$\AA(\O_\pp/P_\pp^s)$ which is induced by an element $\tilde{c}$
of norm one in $\AA(\O_\pp/P_\pp^s)$. Lift it to an element
$\tilde{a}$ of $\AA(\O_\pp)$, then $N(\tilde{a}) \in \O_\pp$ and
is equivalent to $1$ modulo $P_\pp^s$; thus $\tilde{a}$ is
invertible and in $G'(\O_\pp)$.

If $d$ is prime to $\card{\O_\pp/P_\pp}$, then there is some
$\alpha \in \O_\pp$ such that $\alpha^d = \Norm(\tilde{a})$, so
$\Norm(\alpha^{-1}\tilde{a}) = 1$. Thus the automorphism induced
by $\alpha^{-1}\tilde{a}$ on $\AA(\O_\pp)$ is in $G'_1(\O_\pp)$,
call it $a$. Now $a$ covers $c$, and we are done.

When $d$ is not prime to $\card{\O_\pp/P_\pp}$,
we construct an
element of $\AA(\O_\pp)$ by induction. Let $v \in \F_{q^d}$ be an
element with $\tr_{\F_{q^d}/\F_q}(v) = 1$. Take $\tilde{a}_s =
\tilde{a}$. Assume $\Norm(\tilde{a}_i) \equiv 1 \pmod{\pp^i}$ for
some $i \geq s$. Write $\Norm(\tilde{a}_i) \equiv 1 + g \pp^i
\pmod{\pp^{i+1}}$ for some $g \in R$. Let $\tilde{a}_{i+1} = (1-v
g \pp^i)\tilde{a}_i$, so $\tilde{a}_{i+1} \equiv \tilde{a}_i
\pmod{\pp^i}$. Since $\Norm(1 - v g \pp^i) =
\prod_{j=0}^{d-1}(1-\phi^j(vg)\pp^i) \equiv 1- g\pp^i
\pmod{\pp^{i+1}}$ by our choice of $v$, we have
$\Norm(\tilde{a}_{i+1}) \equiv (1-g\pp^i)(1+g\pp^i) \equiv 1
\pmod{\pp^{i+1}}$. Now let $\tilde{a}_0 = \lim {\tilde{a}_i}$.
Then $\Norm(\tilde{a}_0) = \lim {\Norm(\tilde{a}_i)} = 1$, and
$\tilde{a}_0 \equiv \tilde{a}_s \pmod{I}$, proving that the image
of $\tilde{a}_0$ in $\tG'_1(\O_\pp/P_\pp^s)$ is $\tilde{c}$.
Consequently, the image of $\tilde{a}_0$ in $G'_1(\O_\pp)$ covers
$c$.
\end{proof}

Let $H = G'(R,I)$,  $H_1 = H \cap G'_1(R) = G'_1(R,I)$ and
$\Gamma_1 = \Gamma \cap G'_1(R)$.  (Note that the map $\tG'_1(R)
\ra G'(R)$ takes $\tilde{\Gamma}_1 = \tilde{\Gamma} \cap
\tG'_1(R)$ to $\set{\gamma \mul{R}\suchthat \gamma \in
\tilde{\Gamma},\,\Norm(\gamma)=1}$ which may be a proper subgroup
of $\Gamma_1 = \set{\gamma \mul{R}\suchthat \gamma \in
\tilde{\Gamma},\,\Norm(\gamma) \in \mul{R}^d}$).

In order to complete the identification of $\Gamma/\Gamma(I)$, we
need to compute the index of some subgroups of $G'(R)$. {}From the
definition of $R$, it follows that $\mul{R}
=\mul{\F_{q}}\set{y^i(1+y)^j\suchthat i,j \in \Z}$. The norm map
$\mul{\AA(R)} \ra \mul{R}$ is onto, since $\Norm(b_u) = y/(1+y)$,
$\Norm(z) = 1+y$ and the norm $\mul{\F_{q^d}} \ra \mul{\F_{q}}$ is
onto. Therefore, the map $G'(R) \ra \md[d]{R}$, defined by
$a\mul{R}\mapsto \Norm(a){\mul{R}}^d$, is also onto. If $a
\mul{R}$ is in the kernel, then $\Norm(a) = t^d$ for some $t \in
\mul{R}$, and then $\Norm(t^{-1}a) = 1$ and $a\mul{R} = t^{-1} a
\mul{R} \in G'_1(R)$. It follows that $G'(R)/G'_1(R) \isom
\md[d]{R}$, and
\begin{equation}
\dimcol{G'(R)}{G_1'(R)} = \dimcol{\mul{R}}{\mul{R}^d} =
d^2(d,q-1).
\end{equation}
On the other hand, from \Cref{Gamma1},
\Cref{GRGR0} and \Pref{GRT=1}, we obtain
\begin{equation}
\dimcol{G'(R)}{\Gamma_1} = \frac{q^d-1}{q-1}d^2.
\end{equation}

\begin{cor}
$\Gamma_1$ is not contained in any congruence subgroup $G'_1(R,I)$
for $0 \neq I \normali R$ ($I \neq R$).
\end{cor}
\begin{proof}
We may assume $I$ is maximal. By \Tref{isonto} the index of
$G'_1(R,I)$ in $G_1'(R)$ is $\card{\PSL[d](R/I)} \geq
\card{\PSL[d](\F_q)} > \frac{q^d-1}{(q-1)(d,q-1)}$, so that
$\dimcol{G'(R)}{\con{G'_1}{R}{I}} =
d^2(d,q-1)\dimcol{G_1'(R)}{G'_1(R,I)} > \frac{q^d-1}{q-1}d^2
=\dimcol{G'(R)}{\Gamma_1}$.
\end{proof}

\begin{cor}
For every $I \normali R$, $G'_1(R) = \Gamma_1 \cdot H_1$ where
$H_1 = \con{G'_1}{R}{I}$.
\end{cor}
\begin{proof}
By \Pref{GRGR0}, $\Gamma$ and $H$ are normal in $G'(R)$, so
$\Gamma_1 H_1$ is normal in $G'_1(R)$, and $\Gamma_1 H_1 / H_1$ is
normal in $G'_1(R)/H_1 \isom G'_1(R/I)$, which is isomorphic to
$\PSL[d](R/I)$.

Recall that $I = \ideal{p}^s$ where $p$ is prime, and $R/\ideal{p}
= \F_{q^{\alpha}}$ for an appropriate $\alpha$. First assume that
$\PSL[d](\F_{q^\alpha})$ is simple, so that $d
> 2$ or $q^{\alpha} > 3$. The group $\PSL[d](R/I)$ has normal subgroups
$N_j = \Ker(\PSL[d](R/I)\ra \PSL[d](R/\ideal{p}^j))$. Now, the
composition factors are $\PSL[d](R/I)/N_1 \isom
\PSL[d](R/\ideal{p})$, and the $N_{j}/N_{j+1}$ for $j =
1,\dots,s-1$, which are all isomorphic to the zero-trace
subalgebra of $\M[d](\F_{q^{\alpha}})$, and are irreducible
$\PSL[d](\F_{q^\alpha})$-modules. If $N \normal \PSL[d](R/I)$ is a
normal subgroup which is not contained in $N_1$, then
$\PSL[d](R/I) / N_1 \isom N / (N \cap N_1)$ is a group of order
power of $q$, which thus has a normal subgroup of index $q$,
producing an impossible composition factor.

This proves $\Gamma_1 H_1 = G_1$, for otherwise $\Gamma_1 H_1 /
H_1$ is a proper normal subgroup of $G_1/H_1$, which is therefore
contained in $N_1 / H_1 = \con{G'_1}{R}{\sqrt{I}}$, contradicting
the last corollary.

In the special case $d = 2$ and $\card{R/I} = 3$, $G'_1(R)/H_1
\isom A_4$
and has no non-trivial normal subgroups of index $\leq 2 =
\dimcol{G_1'(R)}{\Gamma_1}$. The case $\card{R/I} = 2$ is
impossible since the prime generator $p$ must be different than
$y$ and $1+y$, and so it cannot be linear over $\F_2$.
\end{proof}

There is a delicate point here: the claim that $G'_1(R) = \Gamma_1
\cdot \con{G'_1}{R}{I}$ holds for every proper ideal of $R$, but
cannot be extended to the wider class of ideals of $R_0$ (in the
form $G'_1(R) = \Gamma_1 \cdot \con{G'_1}{R_0}{I}$); indeed, we
saw in the proof of \Tref{main} that $\con{G'_1}{R_0}{I} \sub
\Gamma_1$ for suitable ideals of $R_0$.

{}From this and the definition of $\Gamma_1$ and $H_1$, it follows
that for the subgroups of $G'(R)$ appearing in Figure~\ref{GRsg},
the  intersections and products can be read from the lattice. The
only non-trivial equality is $\Gamma_1 H_1 (\Gamma \cap H) =
\Gamma_1 H \cap \Gamma H_1$, and this too is easily checked (if
$\gamma_1 h = \gamma h_1$ for $\gamma_1 \in \Gamma_1$, $\gamma \in
\Gamma$, $h \in H$ and $h_1 \in H_1$, then $\gamma^{-1} \gamma_1 =
h_1 h^{-1} \in \Gamma \cap H$ so $\gamma_1 (h_1 h^{-1})^{-1} h_1
\in \Gamma_1 H_1 (\Gamma \cap H)$).

In particular, the quotients of the extensions denoted by double
lines are isomorphic to $G_1/H_1 = G'_1(R)/\con{G'_1}{R}{I} \isom
G'_1(R/I) \isom \PSL(R/I)$.

Let $L = R/I$. We now observe that
$$\PSL[d](L) \isom \Gamma_1 H_1 / H_1 \isom \Gamma_1 H / H \sub
\Gamma H / H \sub G'(R)/H \sub \PGL[d](L),$$
so we proved:
\begin{thm}\label{GI}
$\Gamma_I = \Gamma/\Gamma(I) \isom \Gamma H / H$ is a subgroup of
$\PGL(L)$ which contains $\PSL(L)$.
\end{thm}

\ifXY
\begin{figure}
\begin{equation*}
\xymatrix@R=24pt@C=16pt{
    {}
    & {}
    & G'(R) \ar@{-}[d]
    & {}
    & {}
\\
    {}
    & {}
    & \Gamma H \ar@{-}[dr] \ar@{-}[dl]
    & {}
    & {}
\\
    {}
    & \Gamma_1 H \ar@{=}[dl]  \ar@{-}[dr]
    & {}
    & \Gamma H_1 \ar@{-}[dl] \ar@{-}[dr]
    & {}
\\
    H \ar@{-}[dr]
    &
    & \Gamma_1 H \cap \Gamma H_1 \ar@{-}[d] \ar@{=}[dl] \ar@{-}[dr]
    &
    & \Gamma \ar@{-}[dl]
\\
    {}
    & H_1(\Gamma\cap H) \ar@{-}[d] \ar@{-}[dr]
    & G_1' = \Gamma_1 H_1 \ar@{=}[dl]|[l,d]\hole \ar@{-}[dr]|[r,d]\hole
    & \Gamma_1 (\Gamma \cap H) \ar@{-}[d] \ar@{=}[dl]
    &
\\
    {}
    & H_1 \ar@{-}[dr]
    & \Gamma \cap H \ar@{-}[d]
    & \Gamma_1 \ar@{=}[dl]
    & {}
\\
    {}
    & {}
    & \Gamma_1 \cap H_1
    & {}
    & {}
}
\end{equation*}
\caption{Subgroups of $G'(R)$}\label{GRsg}
\end{figure}
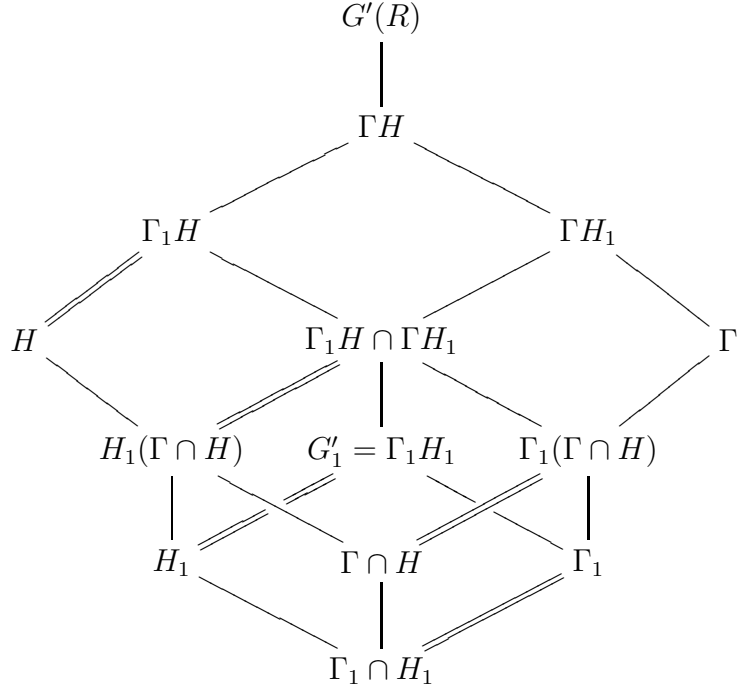
\else
\bigskip
Subgroups of $G'(R)$.
\bigskip
\fi

\medskip

\begin{prop}\label{whatisr}
Let $r$ denote the order of $y/(1+y)$ in $\md{L}$.
Then $\dimcol{\Gamma_I}{\PSL[d](L)} = r$.
\end{prop}
\begin{proof}
The composition $\Gamma_I \isom \Gamma H / H \hookrightarrow
\PGL[d](R/I) \stackrel{\Norm}{\lra} \md[d]{L}$ carries $\Gamma_I$
to the subgroup generated by $\Norm(b) = y/(1+y)$, and its kernel
is $\Gamma_1 H / H \isom \PSL[d](L)$.
\end{proof}

The index of $\PSL[d](L)$ in $\Gamma_I$ can be interpreted in
another way. As mentioned above, the color map $\delta \co \Gamma
\ra \Z/d$ of Figure~\ref{detcol} is in fact an isomorphism
$\Gamma/\Gamma_1 \ra \Z/d$, so it induces an isomorphism $\Gamma
H /\Gamma_1 H \isom \Gamma / \Gamma_1(\Gamma \cap H) \ra \Z/r$. In
other words, there is an $r$-coloring of $\Gamma H / H$.

A remark about colors is in order. For the building $\B$, there is
a $d\times d$ matrix $n_{ij}$ of integers, such that every vertex
of color $i$ has $n_{ij}$ neighbors of color $j$ (in fact $n_{ij}
= \binomq{d}{\abs[]{i-j}}{q}-\delta_{ij}$, where
$\binomq{d}{k}{q}$ is the number of subspaces of dimension $k$ of
$\F_q^d$), and in particular $n_{ii} = 0$ (so the graph is
$d$-partite).

A quotient complex $\B_I$ which has a $d_1$-coloring still
satisfies this property, for the `folded' $d_1\times d_1$ matrix
whose $(i',j')$ entry is the sum of $n_{ij}$ over $i\equiv i'
\pmod{d_1}$ and $j \equiv j' \pmod{d_1}$. Naturally the most
interesting case is when the quotient has $d$ colors, which the
above theorem guarantees quite often:
\begin{cor}
If $y/(1+y)$ has order $d$ in $\md{L}$, then
there is a well defined color
epimorphism $\delta_L \co \Gamma/\Gamma(I) \ra \Z/d\Z$.
\end{cor}

\section{Construction with given index over $\operatorname{PSL}_d(L)$}\label{hitindex}

We combine the results obtained so far to prove Theorem
\ref{mainIntro}, in the following, stronger version (the
assumption $q^e > 4d^2 +1$ is only needed here to control the
index $r$).

\begin{thm}\label{mainIntro+}
Let $q$ be a prime power, $d \geq 2$, $s \geq 1$, $e \geq 1$ ($e >
1$ if $q = 2$), and $r$ a divisor of $(d,q^e-1)$. Assume $d$ is
prime to $q$, or $s = 1$. Also assume that $q^e > 4d^2 + 1$.

Then, for a suitable prime polynomial $\pp \in \F_q[y]$ of degree
$e$ and for $L = \F_q[y]/\ideal{\pp(y)^s}$, the (unique) subgroup
$G$ of $\PGL[d](L)$ which has index $r$ over $\PSL[d](L)$ has a
set of $\binomq{d}{1}{q} + \binomq{d}{2}{q} + \cdots +
\binomq{d}{d-1}{q}$ generators, such that the Cayley graph of $G$
with respect to these generators determines a Ramanujan complex.
\end{thm}

Let $J \normali R$ be a prime ideal and $I = J^s$ for some $s \geq
1$. Let $L_0 = R/J \isom \F_{q^e}$, the residue field of $L =
R/I$.

Notice that $\PGL[d](L)/\PSL[d](L) \isom \md[d]{L}$ is usually not
a cyclic group, so $\Gamma_I$ may not be the full group
$\PGL[d](L)$. However, in two important cases, namely if $d$ is
prime to $q$ or if $L = L_0$, the quotient $\mul{L}/\mul{L}^d
\isom \mul{L_{0}}/\mul{L_{0}}^d$ is cyclic.

We want to choose the generator of $J$ so that $\Gamma/\Gamma(I)$
has a pre-determined index over $\PSL(L)$. Write $J =
\ideal{\pp}$, where
$\pp$ is a prime polynomial different than $y$ and $1+y$, since $R
= \F_q[y,1/y,1/(1+y)]$. In particular, if $q = 2$ then we must
assume $e \geq 2$.

It is convenient to have a uniform bound on the number of non-generators of
a finite field.
\begin{lem}\label{bound}
Let $q$ be a prime power and $e>1$, and let $\F_{q^e}$ be a field
of order $q^e$. There are less than $2q^{e/2}$ elements of
$\F_{q^e}$ which do not generate it over $\F_q$.
\end{lem}
\begin{proof}
We need to bound the size of the union of the maximal subfields of
dimensions $e/p$, for the prime divisors $p$ of $e$. If $e$ is a
prime or equals $4$, then the statement is obvious, so we assume
$e\geq 6$ and is not a prime. Let $C = 9/q$.

We claim that for every prime $p \geq 3$ dividing $e$,
$q^{e/p}<Cp^{-2}q^{e/2}$. Indeed, since $e$ is not a prime, $p\leq
e/2$. Now, in the range $3\leq p \leq e/2$, $h(p) =
2\log(p)-(e/2-e/p)\log(q)$  obtains its maximum on one of the
boundaries, and so is bounded by
$$\max(h(3),h(e/2)) = \max(\log(9)-e\log(q)/6, 2\log(e/2)-(e/2-2)\log(q)).$$
This decreases with $e$, and is thus bounded by its value at
$e=6$, namely $\log(9/q)$.

Since $\sum_{p\geq 3}{p^{-2}}<2/9$ (summing over all the primes), we have that
$\sum_{p|e} q^{e/p} < (1+ C \sum{p^{-2}})q^{e/2} < (1+2/q)q^{e/2}$.
\end{proof}

\begin{prop}\label{Chooseorder}
Let $q,d,s,r$ and $e$ be as in Theorem \ref{mainIntro+}.

If $q^e > 4d^2 + 1$, then there is a polynomial $\pp \in\F_q[y]$
of degree $e$, such that for $J = \ideal{\pp} \normali R$ and $I =
J^s$, $L_0 = R/J$ is the field of order $q^e$, and
$\dimcol{\Gamma_I}{\PSL(L)} = r$ where $L = R/I$.
\end{prop}
\begin{proof}
The assumption that $s = 1$ or $d$ is prime to $q$ guarantees that
$\md[d]{L} \isom \md[d]{L_0}$.

Let $\alpha$ be an element of $\F_{q^e}$ whose order in
$\mul{\F_{q^e}}/\mul{\F_{q^e}}^d$ is $r$. Let $\alpha_1 =
\frac{\alpha}{1-\alpha}$. The map $\alpha\mapsto \alpha_1$ is
one-to-one.

There are $\frac{q^e-1}{(d,q^e-1)}\phi(r) \geq
\frac{q^e-1}{(d,q^e-1)}$ suitable elements $\alpha$, and
$\alpha_1$ belongs to a proper subfield of $\F_{q^e}$ in less than
$2 q^{e/2} \leq \frac{q^e-1}{(d,q^e-1)}$ cases. We can thus assume
$\alpha_1$ generates $\F_{q^e}$ as a field. Then, the minimal
polynomial $f \in \F_q[y]$ of $\alpha_1$ has degree $e$. Take $I =
\ideal{f^s}$. Then $L_0 = R/\ideal{f}$ is isomorphic to
$\F_{q^e}$, and the image of $y$ in $L_0$ satisfies $y/(1+y)
\mapsto \alpha_1/(1+\alpha_1) = \alpha$. It has order $r$ in
$\md{L_0}$, and we are done by \Pref{whatisr}.
\end{proof}

Recall
that the group $\Gamma \sub G'(R)$ (\Dref{Gammadef}) is generated
by the $\binomq{d}{1}{q} = \frac{q^d-1}{q-1}$ elements $b_u =
u(1-z^{-1})u^{-1}$, which can be viewed as $d^2\times d^2$
matrices over $R_0 = \F_q[1/y]$; see \Pref{twodefs}.
The quotient $G'_1(R)/\con{G'_1}{R}{I}$ was computed in
Section~\ref{sec:finq}. For an explicit identification of this
group, it is worth mentioning that $G'(R) / \con{G'}{R}{I} \isom
i_R G'(R) / i_R \con{G'}{R}{I}$ ($i_R$ was defined in Section
\ref{sec2}), and likewise for $G'_1$. Therefore when dealing with
this quotient we may view (the entries of) the generators of
$\Gamma$ as if they were elements of $R/I$.

Let $G$ be the group generated by these matrices over
$R_0/(y^{-e}\pp)^s$. Then $G \isom \Gamma_I$ and, therefore, has
index $r$ over $\PSL[d](L)$. For each $k = 1,\dots,d-1$, the
generators of color $k$ of $G$ are the products $b_{u_1} \ldots
b_{u_k}$, which can be extended to a product of $d$ generators
which equals $1$ (see Section \ref{geom} and \Pref{partofd}); in
particular, every $b_u$ is a generator of color $1$. See
\cite[Sec.~4]{CS1} for an alternative description of the
generators (as elements of $\AA(k)$).

The description of the complex constructed from the Cayley graph
is given at the beginning of Section \ref{sec:finq}. Theorem \ref{main} then
completes the proof of Theorem \ref{mainIntro+}.

\section{Splitting $\AA(R)$}\label{splitx}

So far $\Gamma$ is viewed as a group of $d^2\times d^2$ matrices,
and the representation of $\Gamma$ as $d\times d$ matrices is over
the completion $F$ of $k$. To facilitate the explicit
constructions of finite quotients in the next section, we now
present a splitting extension of dimension $d$ for $G'(k)$ (and
$G'(R)$).

Fix an element $\beta \in \F_{q^d}$ such that $\tr(\beta) \in
\F_q$ is non-zero. Let $\bar{k}_1 = \F_{q^d}(x)$, with the obvious
action of $\Gal(\F_{q^d}/\F_q)$, so the fixed field is $\bar{k} =
\F_q(x)$. Then the norm $\Norm[\F_{q^d}(x)/\F_q(x)](1+\beta x)$ is
a polynomial of degree $d$ in $\F_q[x]$, which has the form $1+t_1
x + \cdots + t_d x^d$ for some $t_i \in \F_q$. Then $t_1 =
\tr(\beta) \neq 0$ by assumption, and $t_d =
\Norm[\F_{q^d}/\F_q](\beta) \neq 0$, since $\beta \neq 0$.

We embed $k = \F_q(y)$ into $\bar{k}$ by mapping
\begin{equation} \label{ydef}
y \mapsto t_1 x + \cdots + t_d x^d.
\end{equation}
Now $k$ becomes a
subfield of codimension $d$ of $\F_q(x)$.

\begin{prop}\label{Fxy}
Taking the completion, $\F_q\db{x} = \Fy$, and even $\F_q[[x]] =
\F_q[[y]]$.
\end{prop}
\begin{proof}
Before giving the proof, we should explain the notation. By
$\F_q\db{x}$ we mean, as usual, the field of Taylor series
$\set{\sum_{i=-v}^{\infty}{\alpha_i x^i}}$ with coefficients from
$\F_q$; $\F_q[[x]]$ is the subring of elements of the form
$\sum_{i=0}^{\infty}{\alpha_i x^i}$. As abstract fields, it
trivially holds that $\Fy \isom \F_q\db{x}$. But now $\Fy$ is
the subfield of $\F_q\db{x}$ consisting of Taylor series in the
element $y$ of \Eq{ydef}.

Working in $\Fy$,  \Eq{ydef} has a solution for $x$
by Hensel's lemma (since the derivation at zero, namely $t_1$, is
non-zero). In fact, the solution has the form $x = \frac{1}{t_1}y
- \frac{t_2}{t_1^3}y^2 + \dots$, so that $x \in \F_q[[y]]$.
\end{proof}
Tensoring with $\F_{q^d}$, we also have that $\Fy[q^d] =
\F_{q^d}\db{x}$. Define $\bar{R}_T = \F_q[x,1/(1+y)]$ and $\bar{R}
= \F_q[x,1/y,1/(1+y)]$, so that $R_T \sub \bar{R}_T$ and $R \sub
\bar{R}$. Since $1+y$ is invertible in $\F_q[[y]]$, $\bar{R}_T
\sub \F_q[[y]]$. The rings mentioned so far appear in
Figure~\ref{fi:fields}.

\ifXY
\begin{figure}
\begin{equation*}
\xymatrix@R=16pt@C=6pt{
    {}
    & {}
    & {}
    & {}
    & {}
    & {}
\\
    {}
    & {}
    & {}
    & {\Fy[q^d]} \ar@{-}[dr] \ar@{--}[dl]
    & {}
    & {}
\\
    {}
    & {}
    & {\F_{q^d}(x)} \ar@{-}[dr] \ar@{-}[dl]
    & {}
    & {F = \Fy} \ar@{--}[dl] \ar@{-}[dd]
    & {}
\\
    {}
    & {k_1 = \F_{q^d}(y)} \ar@{-}[dr]
    & {}
    & {\F_{q}(x)} \ar@{-}[dl] \ar@{-}[d]
    & {}
    & {}
\\
    {}
    & {}
    & {k = \F_{q}(y)} \ar@{-}[d]
    & {\bar{R}} \ar@{-}[dl] \ar@{-}[d]
    & {\O = \F_q[[y]]} \ar@{--}[dl]
    & {}
\\
    {}
    & {}
    & {R} \ar@{-}[d]
    & {\bar{R}_T}
    & {}
    & {}
\\
    {}
    & {}
    & {R_T} \ar@{-}[ur]
    & {}
    & {}
    & {}
}
\end{equation*}
\caption{Subrings of $\F_{q^d}\db{y}$}\label{fi:fields}
\end{figure}
\else
\bigskip
Zoo.
\bigskip
\fi

Notice that $\F_q(x)$ is a splitting field of $\AA(k)$, as $1+y =
\Norm(1+\beta x)$ is a norm in the extension
$\F_{q^d}(x)/\F_q(x)$. We thus obtain the chain of
inclusions
\begin{eqnarray}\label{inc}
\AA(R) & \sub & \AA(k) \,\hra\, \AA(k) \tensor[k] \F_q(x) \\
& \isom & \M(\F_q(x)) \,\hra\, \M(\F_q\db{x}) = \M(\Fy), \nonumber
\end{eqnarray}
refining \Eq{ARinMF}. The
last equality in this chain follows from \Pref{Fxy}.

In fact, $\AA(\bar{R}_T)$ is already split,
and we now present an explicit isomorphism $\AA(\bar{R}_T) \isom
\M(\bar{R}_T)$.

\begin{prop}\label{gse}
Let $M = \bar{R}_T\tensor \F_{q^d}$, and let $\rho \co M \ra
\End_{\bar{R}_T}(M)$ be the regular representation, defined by
$\rho \co r\tensor u \mapsto \rho_{r} \tensor \rho_{u}$, where
$\rho_{r}(r') = rr'$ for $r' \in \bar{R}_T$, and $\rho_{u}(u') =
uu'$ for $u \in \F_{q^d}$.

The map $\rho$ is extended by $z \mapsto 1\tensor \phi + \rho_x \tensor
\rho_{\beta} \phi$ to an isomorphism
$$\AA(\bar{R}_T) \stackrel{\sim\,}{\lra}
\End_{\bar{R}_T}(M).$$
\end{prop}
\begin{proof}
By definition $\rho(z) = \rho(1+\beta x) \cdot (1\tensor \phi)$,
so that $\rho(z)^{-1} = (1\tensor \phi^{-1})\cdot \rho(1+\beta
x)^{-1}$. We need to check the defining relations of
$\AA(\bar{R}_T)$. Let $r \tensor u, r' \tensor u' \in M$; then
\begin{eqnarray*}
\rho(z)\rho(r \tensor u) (r' \tensor u') & = &
    \rho(z)(r r' \tensor u u') \\
& = & (1+\beta x) (r r' \tensor \phi(u u')) \\
& = & (1+\beta x) (r \tensor \phi(u))  (r'\tensor \phi(u')) \\
& = & \rho(r \tensor \phi(u)) (1+\beta x) (r'\tensor \phi(u')) \\
& = & \rho(r \tensor \phi(u)) \rho(z) (r' \tensor u'),
\end{eqnarray*}
so that $\rho(z) \rho(r \tensor u)\rho(z)^{-1} = \rho(r \tensor
\phi(u))$. Also, $\rho(z)^d = (\rho(1+\beta x) \cdot (1\tensor
\phi))^d = \rho(\Norm(1+\beta x)) \cdot (1\tensor \phi)^d =
\rho(1+y)$.

To check that the extended $\rho$ is an embedding, let $\bar{k} =
\F_q(x)$, and recall that $\AA(\bar{k}) \isom
\End_{\bar{k}}(\bar{k} \tensor \F_{q^d})$ by the norm condition.
The result then follows from the commutativity of the diagram in
Figure \ref{fi:di}, and the fact that $\AA(\bar{R}_T)$ and
$\End_{\bar{R}_T}(M)$ have the same rank
$d^2$ over $\bar{R}_T$.

\ifXY
\begin{figure}
\begin{equation*}
\xymatrix{
    {\AA(\bar{k})} \ar@{->}[r]^(0.35){\isom}
    & {\End_{\bar{k}}(\bar{k}\tensor \F_{q^d})}
\\
    {\AA(\bar{R}_T)} \ar@{^(->}[u] \ar[r]
    & {\End_{\bar{R}_T}(\bar{R}_T\tensor \F_{q^d})} \ar@{^(->}[u]
}
\end{equation*}
\caption{Maps in \Pref{gse}}\label{fi:di}
\end{figure}
\else
\bigskip
Zoo 2.
\bigskip
\fi

\end{proof}

Since $\bar{R}_T \tensor \F_{q^d}$ is a free module of rank $d$
over $\bar{R}_T$, we obtain an isomorphism  $\AA(\bar{R}_T) \isom
\M(\bar{R}_T)$.

\begin{exmpl}\label{gse-exmp}
Suppose $q = 7$ and $d = 3$, and let $\alpha = \sqrt[3]{2}$ over
$\F_7$, so that $\F_{7^3} = \F_7[\alpha]$. If $\phi$ is the
Frobenius automorphism, we have that $\phi(\alpha) = \alpha^7 =
4\alpha$. As a normal basis of $\F_{7^3}/\F_7$ we take $\zeta_0 =
1+\alpha+\alpha^2$ and $\zeta_i = \phi^{i}(\zeta_0)$. We choose
$\beta = \alpha - 2$ (so that $\tr(\beta) = 1$). Of course
$\set{\zeta_i}$ form a basis of $\bar{R}_T\tensor \F_{7^3}$ as a
free module over $\bar{R}_T$. Computing directly, we find that
$\Norm(1+\beta x) = 1+x-2x^2+x^3$, so that $y = x-2x^2+x^3$.

In order to present an embedding of $\AA(\bar{R}_T) = \F_{7^3}[x,1/(1+y)][z]$
into $\M[3](\F_{7}[x,1/(1+y)])$, we need to describe the images of
$\F_{7^3}$, $x$, $1/(1+y)$ and $z$. For the subfield, check that
$$
\alpha \mapsto \(\begin{matrix} 6 & 3 & 6 \\ 5 & 5 & 6 \\ 5 & 3 & 3 \end{matrix}\)
$$
defines the embedding $\F_{7^3} \ra \M[3](\F_7)$ via the basis we
chose.
\def\pl{{+}}
Secondly, $x$ and $1/(1+y)$ are scalars in $\bar{R}_T$ and so are mapped to
scalar matrices. Finally, since $\phi(\zeta_i) = \zeta_{i+1}$, by
the above proposition we have
\begin{eqnarray*}
z & \mapsto &\rho(1+\beta x) \cdot (1\tensor \phi) \\
& = & \(\begin{matrix} 1\pl 4x & 3x & 6x \\
5x & 1\pl 3x & 6x
\\ 5x & 3x & 1\pl x
\end{matrix}\)
\cdot
\(\begin{matrix} 0 & 0 & 1 \\ 1 & 0 & 0 \\ 0 & 1 & 0 \end{matrix}\)
=
\(\begin{matrix} 3x & 6x & 1\pl 4x\\ 1\pl 3x & 6x & 5x\\  3x &
1\pl x & 5x
\end{matrix}\)
.
\end{eqnarray*}
One can check that over $\bar{R}_T$, $z^3$ is the scalar matrix
$1+y$.
\end{exmpl}

\section{Explicit constructions}\label{sec:explicit}

In this section we show how to choose the generator $\pp^s \in \F_q[y]$
of $I$, such that $\Gamma_I = \Gamma/\Gamma(I)$ will have an explicit
representation as $d\times d$ matrices over $L$. The construction
works for arbitrary $q$, $d$, $e = \deg(\pp)$ and $s \geq 1$, but we
do not try to control the index $r$ of $\Gamma_I$ over $\PSL[d](L)$
(which can be computed by \Pref{whatisr}).

Recall that $H = \con{G'}{R}{I}$. Our current realization of $\Gamma H / H$ as a subgroup of $\PGL[d](R/I)$
 is inexplicit, as we use the abstract
isomorphism $G'(R/I) \isom \PGL(R/I)$. We do have an explicit embedding of
$i_R$ of $G'(R)$ into $d^2\times d^2$ matrices over $R$, as described before
\Pref{twodefs}, but for practical reasons it is more convenient to have
a representation of dimension $d$.

In Proposition \ref{gse} (and Example \ref{gse-exmp}) we gave an
explicit isomorphism of $\AA(\bar{R}_T)$ to matrices over
$\bar{R}_T$, which in particular embeds $G'(R) \hra \PGL[d](\bar{R})$. Given $I \normali R$, there is an ideal $\bar{I} \normali \bar{R}$
such that $I = \bar{I} \cap R$ (since $\bar{R} = R[x]$
is an integral extension over $R$ and $\bar{R}$ is a  principal ideal domain).
Let $\bar{L} = \bar{R}/\bar{I}$, then there is a natural embedding $L \sub
\bar{L}$, and we thus have an embedding $G'(L) \hra \PGL(\bar{L})$.

We will show that it is possible to choose $I$ so that
$\bar{L} \isom L$, where the image of $x$
in $L$ (which is what we need for the splitting map of \Pref{gse}) is
explicitly identified.

We need some details on the possible choices of $\bar{I}$. Write
$y(\lam) = t_1 \lam+ \cdots + t_d \lam^d$ for the polynomial
defined by \Eq{ydef}, so $y(x)$ is our usual element $y$ of
$\F_q[x]$. If $\pp(\lam) \in \F_q[\lam]$ is irreducible and not
equal to $\lam$ or $\lam+1$, then $\pp(y) \in R$ generates a prime
ideal $J \normali R$. Since the embedding of $R$ into $\bar{R}$ is
by sending $y$ to $y(x)$, $J$ lifts to the ideal $\bar{J} =
\ideal{\pp(y(x))}$ of $\bar{R}$. Let $\pp(y(\lam)) = g_1(\lam)
\ldots g_t(\lam)$ be the decomposition into prime factors over
$\F_q$ and $\bar{J}_i = \ideal{g_i(x)} \normali \bar{R}$, then
$\ideal{\pp(y(x))} = \bar{J}_1 \ldots \bar{J}_t$, where each
$\bar{J}_i$ is a prime ideal of $\bar{R}$ covering $J$. If
$\alpha$ is a root of $g_i$ in some extension, then $y(\alpha) \in
\F_q[\alpha]$ is a root of $\pp(\lam)$. Therefore $\pp$ splits in
$\F_q[\alpha]$ and thus $\deg(\pp(\lam)) \divides
\dimcol{\F_q[\alpha]}{\F_q} = \deg(g_i)$.  This shows once more
that $\bar{L}^i_0 = \bar{R}/\bar{J}_i$ is an extension of $L_0 =
R/J$, and in fact $\dimcol{\bar{L}^i_0}{L_0}=
\deg(g_i)/\deg(\pp)$. We remark that the degrees of the $g_i$ are
not necessarily equal.

Let $\gamma_j$ be the roots of $\pp(\lam)$ in $L_0$, then $\pp(\lam) = \prod(\lam-\gamma_j)$ and $\pp(y(\lam)) = \prod(y(\lam) - \gamma_j)$. In particular the roots
of $\pp(y(\lam))$ are roots of the factors $y(\lam) - \gamma_j$.
Let $\alpha$ be a root of $\pp(y(\lam))$ coming from an irreducible factor
$g_i(\lam)$. Then $\bar{L}^i_0 = \F_q[\alpha] = L_0[\alpha]$ (since a generator
$\gamma_j$ of $L_0$ equals $y(\alpha)$). It follows that $\pp(y(\lam))$ has a
factor $g_i$ of degree $\deg(\pp)$ iff $\bar{L}_0^i = L_0$,
iff $y(\lam) - \gamma_j$ has a root in $L_0$
for some $j$ (a condition which is easily seen to be independent of $j$).

We now show that $J = \ideal{\pp(y)} \normali R$ can be chosen so
that $\bar{L}_0 = L_0$ and $\bar{L} = L$, which provides an
explicit realization of \Tref{GI}.

\begin{thm}\label{LbarL}
Let $q$ be a prime power and $d$ an integer. Assume
$q^e \geq 4d^2$.
Then there are irreducible polynomials $\pp,g \in \F_q[\lam]$ of degree $e$
such that for every $s \geq 1$, the embedding $R \ra \bar{R}$
induces an isomorphism $R/I \isom \bar{R}/\bar{I}$ for
$I = \ideal{\pp(y)^s} \normali R$ and
$\bar{I} = \ideal{g(x)^s} \normali \bar{R}$.
\end{thm}
\begin{proof}
Let $\F_{q^e}$ be the field of order $q^e$. For every $\alpha \in
\F_{q^e}$, consider the element $y(\alpha)$. Since $\deg(y(\lam))
= d$, there are at most $d$ elements $\alpha$ with the same
$y(\alpha)$, so in particular (by Lemma \ref{bound}), $y(\alpha)$
is not a generator of $\F_{q^e}$ for less than  $2d q^{e/2}$
choices of $\alpha$. By the assumption on $e$, we can choose
$\alpha$ for which $\gamma = y(\alpha)$ generates the field of
$q^e$ elements (so necessarily $\alpha$ is also a generator). Let
$\pp(\lam)$ and $g(\lam)$ be the minimal polynomials of $\gamma$
and $\alpha$ over $\F_q$, respectively. Since $\alpha$ and
$\gamma$ are generators, $\deg(\pp) = \deg(g) = e$.

Now let $I = \ideal{\pp(y)^s}\normali R$ and $\bar{I} =
\ideal{g(x)^s} \normali \bar{R}$. Since $g(x) \divides \pp(y(x))$
and $g(x)$ is prime, $I = \bar{I} \cap R$, and so $L = R/I$ embeds
into $\bar{L} = \bar{R}/\bar{I}$. Since $\deg(\pp) = \deg(g)$,
$\card{L} = \card{\bar{L}}$, and so they are isomorphic.
\end{proof}

The same choice of $\pp$ and $g$ works for every $s \geq 1$, so
the proposition shows that in fact the fraction field $\bar{k} =
k[x]$ of $\bar{R}$ is contained in the completion $k_\pp$.

We conclude this section with a detailed algorithm giving the
generators for a Cayley complex of $\PGL[d](L)$. Let $q$ be a
given prime power, $d \geq 2$, and $e \geq 1$. Assume $q^e \geq
4d^2$.

\begin{algo}
{\rm \begin{enumerate} \item In practice, the elements of
$\F_{q^d}$ are polynomials of degree at most $d$ over $\F_q$,
modulo a fixed irreducible polynomial of degree $d$. Fix a basis
$\zeta_0,\dots,\zeta_{d-1}$ for $\F_{q^d}$ over $\F_q$.
Let $\phi$ be any generator of $\Gal(\F_{q^d}/\F_q)$ (namely
exponentiation by $q^{\ell}$ where $\ell$ is any fixed integer prime to
$d$).
Let $\beta \in \F_{q^d}$ be an element with
$\beta+\dots+\phi^{d-1}(\beta) \neq 0$, let $t_1,\dots,t_d \in
\F_q$ be defined by $(1+\beta \lam)\dots (1+\phi^{d-1}(\beta)
\lam) = 1+ t_1 \lam+\cdots + t_{d} \lam^{d}$, and set $y(\lam) =
t_1 \lam+ \cdots + t_{d} \lam^{d}$. For an element $c \in
\mul{\F_{q^d}}$, $\rho(c) \in \GL[d](\F_q)$ is the matrix defined
by $c \cdot \zeta_j = \sum_{i=0}^{d-1}{(\rho(c))_{ij} \zeta_i}$.
Likewise, let $\varphi \in \GL[d](\F_q)$ be the matrix defined by
$\phi(\zeta_j) = \sum_{i=0}^{d-1}{\varphi_{ij} \zeta_i}$.

\item Let $\alpha \in \F_{q^e}$ be an element such that $\gamma =
y(\alpha)$ generates $\F_{q^e}$ as a field (the existence of which
is guaranteed by \Tref{LbarL}). A quick way to check this property
is to verify that $\gamma^{q^{e'}} \neq \gamma$ for every proper
divisor $e'$ of $e$. Let $p(\lam)$ and $g(\lam)$ be the minimal
polynomials (of degree $e$) of $\gamma$ and $\alpha$,
respectively, over $\F_q$ (found by Gaussian elimination on
$1,\gamma,\dots,\gamma^{e}$ or $1,\alpha,\dots,\alpha^{e}$).

\item Fix $s \geq 1$. Let $\bar{L}$ denote the ring of polynomials of
degree $< es$ in the variable $x$, with operations modulo
$g(x)^s$. Let $\rho(z) \in \GL[d](L)$ be the matrix
$(1+\rho(\beta)x)\cdot \varphi$, as in \Pref{gse} and the example
following it. Let $b = 1-\rho(z)^{-1}$ (Gauss elimination can be
used to invert $\rho(z)$, where the inverse of $f(x) \in L$ is
computed via the extended Euclid's algorithm for $f(\lam)$ and
$g(\lam)^s$). For every $u \in \mul{\F_{q^d}}/\mul{\F_q}$ (one
representative from each class), set $b_u = \rho(u)b\rho(u)^{-1}$.

\item The matrices $\set{b_u}$ generate a subgroup $G$ of
$\PGL[d](L)$, which contains $\PSL[d](L)$. The index $r =
\dimcol{G}{\PSL[d](L)}$ can be computed using \Pref{whatisr}.

\item Let $P = \set{(u_1,\dots,u_d) \suchthat b_{u_1} \dots b_{u_d} = 1}$.
For every $k = 1,\dots, d-1$, let $S_k$ denote the set of products
$b_{u_1}\ldots b_{u_k}$ for the headers $(u_1,\dots,u_k)$ of
vectors in $P$. There will be $\binomq{d}{k}{q}$ different
products (up to scalar multiples over $L$), and $S_{d-k}$ is the
set of inverses of $S_k$.

\item The Cayley complex of $G$ with respect to $S_1 \cup \ldots \cup
S_{d-1}$ defines a simplicial complex (the clique complex of the
Cayley graph). This complex is Ramanujan by \Tref{main}.
\end{enumerate}
}
\end{algo}

\section{An example}\label{sec:example}

We demonstrate the construction of the group $\Gamma$, and the
embedding into $d^2\times d^2$ and $d \times d$ matrices. We
choose $q = 2$ and $d = 3$, and write $\F_8 = \F_2[v \subjectto
v^3=v+1]$.  Fix the basis $\set{\zeta_0,\zeta_1,\zeta_2} =
\set{1,v,v^2}$ for $\F_8$ over $\F_2$. Then $\set{\zeta_i z^j}$ is
a basis for $\AA(R) = R[v,z \subjectto v^3=v+1,\,z^3=1+y,\,zvz^{-1}
= v^2]$ over $R =\F_2[y,1/y,1/(1+y)]$. By definition, $\Gamma \sub
\mul{\AA(R)}$ is generated by the conjugates $u(1-z^{-1})u^{-1}$
for $u \in \mul{\F_{8}}/\mul{\F_{2}} = \mul{\F_{8}} = \sg{v}$. These
generators act on $\AA(R)$ by conjugation, and the corresponding
$9\times 9$ matrices over $R$ are given in the right-hand column
of Table~\ref{ExampleTab} (they are all of the form $a+by^{-1}$
for $a,b\in \GL[9](\F_2)$).

The relations of $\Gamma$ were discussed in Section \ref{geom}. If
we set $b_i = v^i(1-z^{-1})v^{-i}$, conjugation by $v$ induces an
outer automorphism of order $7$ (namely $b_i \mapsto b_{i+1}$),
and the defining relations are
$$b_0 b_3 = b_4 b_2 = b_6 b_5,$$
$$b_0 b_5 = b_2 b_1 = b_3 b_6,$$
$$b_0 b_6 = b_1 b_4 = b_5 b_3,$$
$$b_0 b_3 b_1 = 1,$$
and their $v$-conjugates.

We continue to find a representation of dimension $3$ of a finite
quotient of $\Gamma$. The first step is to define $R$ as a subring
of $\bar{R} = R[x]$, so we need to choose $\beta$ with $\tr(\beta)
\neq 0$; here we take $\beta = 1+v$, so that $1+y =
\Norm[\F_8/\F_2](1+\beta x) = 1 + x + x^3$ and $y(x) = x+x^3$.
Then $R \hra \bar{R}$ by sending $y \mapsto y(x)$.

Mimicking Example \ref{gse-exmp}, the embedding $\AA(R) \hra
\PGL[3](\bar{R})$ is  defined by
$$
v \mapsto \(\begin{matrix} 0 & 0 & 1 \\ 1 & 0 & 1 \\ 0 & 1 & 0 \end{matrix}\)
$$
and
\begin{eqnarray*}
z & \mapsto &\rho(1+\beta x) \cdot (1\tensor \phi) \\
& = & \(\begin{matrix} 1{{+}}x & 0 & x \\
x & 1{{+}} x & x
\\ 0 & x & 1{{+}}x
\end{matrix}\)
\cdot
\(\begin{matrix} 1 & 0 & 0 \\ 0 & 0 & 1 \\ 0 & 1 & 1 \end{matrix}\)
=
\(\begin{matrix} 1{{+}}x & x & x\\ x & x & 1\\
0 & 1{{+}} x & 1
\end{matrix}\)
.
\end{eqnarray*}

Check that $z^3$ is mapped to the scalar matrix $1+x+x^3 =
1+y(x)$, so $z^{-1} = \frac{1}{1+y(x)}z^2$, and $b = 1-z^{-1}
\mapsto
\frac{1}{1+y(x)}\(\begin{matrix} x+x^3 & x^2 & x+x^2 \\ x & x^3 & 1+x+x^2 \\
x+x^2 & 1+x^2 & 1+x^3 \end{matrix}\)$, where $1+y(x)$ is
invertible in $\bar{R}$ by definition. We emphasize that this
embedding is not onto (since $\AA(R)$ does not split).

Finally, we present a finite quotient of $\Gamma$. Choose $e = 4$,
so $\Gamma$ will map onto $\PGL[3](L)$, where $L$ is a finite
local ring whose residue field is $\F_{16}$ (in general $\Gamma$
could be mapped onto a subgroup of $\PGL[3](L)$ containing
$\PSL[3](L)$). Applying the proof of \Tref{LbarL} and the
algorithm given in Section \ref{sec:explicit}, we look for an
element $\alpha \in \F_{16}$ such that $y(\alpha)$ is a generator.

Arbitrarily choosing an irreducible polynomial of degree $4$ over
$\F_2$, we set $\F_{16} = \F_2[t \subjectto t^4+t+1=0]$. It then
turns out that $y(t) = t^3+t$ is a generator (since $y(t^2) \neq
y(t)$), so we choose $\alpha = t$ and $\gamma = y(t)$. The minimal
polynomial of $\alpha$ is of course $g(\lam) = \lam^4+\lam+1$, and
the minimal polynomial of $\gamma$ is computed to be $\pp(\lam) =
\lam^4+\lam^3+\lam^2+\lam+1$. Now $g(\lam)$ divides $\pp(y(\lam))
= \pp(\lam^3+\lam) =
\lam^{12}+\lam^9+\lam^7+\lam^6+\lam^5+\lam^4+\lam^2+\lam+1 =
g(\lam)(\lam^8+\lam^4+\lam^3+\lam^2+1)$. For any $s \geq 1$, when
we take $I = \ideal{\pp(y)^s}$ and $\bar{I} = \ideal{g(x)^s}$, we
have $R \cap \bar{I} = I$. Since $\deg(\pp) = \deg(g)$, $R/I \isom
\bar{R}/\bar{I} = \F_2[x]/\sg{(x^4+x+1)^s}$.

The index over $\PSL$ is the order of $y/(1+y) = y^3+y+1$ in
$\F_2[y]/\ideal{y^4+y^3+y^2+y+1}$ modulo its cubes; the order is
$3$, so $\Gamma$ maps onto $\PGL(L)$. The $3\times 3$ matrices
over $\bar{R}$ corresponding to the generators $b_i$ of $\Gamma / (\con{G'}{R}{I} \cap \Gamma) \isom \PGL[3](R/I)$
are given in the intermediate column of Table~\ref{ExampleTab}. The Cayley
graph with respect to these generators is Ramanujan by Theorem \ref{main}.

\newcommand\doMATnine[2]{{\begin{tiny}${#1} + \frac{1}{y} {#2}$\end{tiny}}}
\newcommand\MATnine[9]{{\(\!\!\! \begin{array}{ccccccccc}#1 \\[-0.14cm] #2 \\[-0.14cm] #3 \\[-0.04cm] #4 \\[-0.14cm] #5 \\[-0.14cm] #6\\[-0.04cm] #7\\[-0.14cm] #8\\[-0.14cm] #9\end{array}\!\!\! \)\! }}
\newcommand\Elemnine[1]{{\!\!\!\!\,{#1}\!\!\;\!\!}}
\newcommand\Linenine[9]{{#1}\!\! &\Elemnine{#2}&\Elemnine{#3}\,&\Elemnine{#4}&\Elemnine{#5}&\Elemnine{#6}\,&\Elemnine{#7}&\Elemnine{#8} & \!\! {#9}}
\newcommand\MATthree[3]{{\(\!\!\! \begin{array}{ccc}#1 \\[-0.0cm] #2 \\[-0.0cm] #3\end{array}\!\!\! \)\!}}
\newcommand\Linethree[3]{{#1}\! &\!\! {#2}\! &\!\! {#3}}

\begin{figure}
\begin{center}
\begin{tabular}{l|c|c}
{} & $b_u \in \PGL[3](\bar{R})$ & $b_u \in \GL[9](R)$
\\
 \hline
$b_0$ & $\MATthree{\Linethree{x{+}x^3}{x^2}{x{+}x^2}}{
\Linethree{x}{x^3}{1{+}x{+}x^2}}{
\Linethree{x{+}x^2}{1{+}x^2}{1{+}x^3}}$ &
\doMATnine{
\MATnine{ \Linenine{1}{0}{0}{0}{0}{0}{0}{0}{0}}{
\Linenine{0}{1}{0}{0}{0}{1}{0}{0}{1}}{
\Linenine{0}{0}{1}{0}{1}{1}{0}{1}{1}}{
\Linenine{0}{0}{0}{1}{0}{0}{0}{0}{0}}{
\Linenine{0}{0}{0}{0}{1}{0}{0}{0}{1}}{
\Linenine{0}{0}{0}{0}{0}{1}{0}{1}{1}}{
\Linenine{0}{0}{0}{0}{0}{0}{1}{0}{0}}{
\Linenine{0}{0}{0}{0}{0}{0}{0}{1}{0}}{
\Linenine{0}{0}{0}{0}{0}{0}{0}{0}{1} }}{\MATnine{
\Linenine{0}{0}{0}{0}{0}{0}{0}{0}{0}}{
\Linenine{0}{0}{1}{0}{0}{1}{0}{0}{1}}{
\Linenine{0}{1}{1}{0}{1}{1}{0}{1}{1}}{
\Linenine{0}{0}{0}{0}{0}{0}{0}{0}{0}}{
\Linenine{0}{0}{1}{0}{0}{1}{0}{0}{1}}{
\Linenine{0}{1}{1}{0}{1}{1}{0}{1}{1}}{
\Linenine{0}{0}{0}{0}{0}{0}{0}{0}{0}}{
\Linenine{0}{0}{1}{0}{0}{1}{0}{0}{1}}{
\Linenine{0}{1}{1}{0}{1}{1}{0}{1}{1} }}
\\
$b_1$ &
$\MATthree{\Linethree{1{+}x{+}x^2{+}x^3}{x{+}x^2}{1{+}x^2}}{\Linethree{1{+}x}{x^2{+}x^3}{1}}
{\Linethree{1{+}x^2}{x}{x^3}}$ &
\doMATnine{
\MATnine{ \Linenine{1}{0}{0}{0}{0}{0}{0}{0}{0}}{
\Linenine{0}{1}{0}{0}{1}{0}{0}{1}{1}}{
\Linenine{0}{0}{1}{1}{1}{1}{1}{0}{0}}{
\Linenine{0}{0}{0}{1}{0}{0}{0}{1}{1}}{
\Linenine{0}{0}{0}{0}{1}{0}{1}{0}{0}}{
\Linenine{0}{0}{0}{0}{0}{1}{0}{0}{0}}{
\Linenine{0}{0}{0}{0}{0}{0}{1}{0}{0}}{
\Linenine{0}{0}{0}{0}{0}{0}{0}{1}{0}}{
\Linenine{0}{0}{0}{0}{0}{0}{0}{0}{1} }}{\MATnine{
\Linenine{0}{0}{0}{0}{0}{0}{0}{0}{0}}{
\Linenine{0}{0}{1}{0}{1}{0}{0}{1}{1}}{
\Linenine{0}{1}{1}{1}{1}{1}{1}{0}{0}}{
\Linenine{0}{0}{1}{0}{1}{0}{0}{1}{1}}{
\Linenine{0}{1}{1}{1}{1}{1}{1}{0}{0}}{
\Linenine{0}{0}{0}{0}{0}{0}{0}{0}{0}}{
\Linenine{0}{1}{0}{1}{0}{1}{1}{1}{1}}{
\Linenine{0}{0}{1}{0}{1}{0}{0}{1}{1}}{
\Linenine{0}{1}{0}{1}{0}{1}{1}{1}{1} }}
\\
$b_2$ & $\MATthree{ \Linethree{1{+}x^2{+}x^3}{1{+}x^2}{x} }{
\Linethree{1{+}x{+}x^2}{x{+}x^3}{x^2} }{
\Linethree{x}{1{+}x}{x^2{+}x^3} }$ &
\doMATnine{
\MATnine{ \Linenine{1}{0}{0}{0}{0}{0}{0}{0}{0}}{
\Linenine{0}{1}{0}{1}{0}{1}{1}{0}{0}}{
\Linenine{0}{0}{1}{1}{1}{0}{1}{0}{1}}{
\Linenine{0}{0}{0}{1}{0}{0}{0}{0}{1}}{
\Linenine{0}{0}{0}{0}{1}{0}{0}{0}{0}}{
\Linenine{0}{0}{0}{0}{0}{1}{1}{0}{0}}{
\Linenine{0}{0}{0}{0}{0}{0}{1}{0}{0}}{
\Linenine{0}{0}{0}{0}{0}{0}{0}{1}{0}}{
\Linenine{0}{0}{0}{0}{0}{0}{0}{0}{1} }}{\MATnine{
\Linenine{0}{0}{0}{0}{0}{0}{0}{0}{0}}{
\Linenine{0}{0}{1}{1}{0}{1}{1}{0}{0}}{
\Linenine{0}{1}{1}{1}{1}{0}{1}{0}{1}}{
\Linenine{0}{1}{0}{0}{1}{1}{0}{0}{1}}{
\Linenine{0}{0}{0}{0}{0}{0}{0}{0}{0}}{
\Linenine{0}{0}{1}{1}{0}{1}{1}{0}{0}}{
\Linenine{0}{1}{1}{1}{1}{0}{1}{0}{1}}{
\Linenine{0}{1}{1}{1}{1}{0}{1}{0}{1}}{
\Linenine{0}{0}{1}{1}{0}{1}{1}{0}{0} }}
\\
$b_3$ &
$\MATthree{\Linethree{x{+}x^2{+}x^3}{x}{1{+}x}}{\Linethree{1}{1{+}x{+}x^2{+}x^3}{x{+}x^2}}
{\Linethree{1{+}x}{1{+}x{+}x^2}{x{+}x^3}}$ &
\doMATnine{
\MATnine{ \Linenine{1}{0}{0}{0}{0}{0}{0}{0}{0}}{
\Linenine{0}{1}{0}{0}{1}{1}{1}{0}{1}}{
\Linenine{0}{0}{1}{1}{0}{0}{1}{1}{0}}{
\Linenine{0}{0}{0}{1}{0}{0}{0}{1}{1}}{
\Linenine{0}{0}{0}{0}{1}{0}{1}{0}{1}}{
\Linenine{0}{0}{0}{0}{0}{1}{0}{1}{1}}{
\Linenine{0}{0}{0}{0}{0}{0}{1}{0}{0}}{
\Linenine{0}{0}{0}{0}{0}{0}{0}{1}{0}}{
\Linenine{0}{0}{0}{0}{0}{0}{0}{0}{1} }}{\MATnine{
\Linenine{0}{0}{0}{0}{0}{0}{0}{0}{0}}{
\Linenine{0}{0}{1}{0}{1}{1}{1}{0}{1}}{
\Linenine{0}{1}{1}{1}{0}{0}{1}{1}{0}}{
\Linenine{0}{1}{0}{1}{1}{1}{0}{1}{1}}{
\Linenine{0}{0}{1}{0}{1}{1}{1}{0}{1}}{
\Linenine{0}{1}{0}{1}{1}{1}{0}{1}{1}}{
\Linenine{0}{1}{0}{1}{1}{1}{0}{1}{1}}{
\Linenine{0}{0}{0}{0}{0}{0}{0}{0}{0}}{
\Linenine{0}{0}{1}{0}{1}{1}{1}{0}{1} }}
\\
$b_4$ &
$\MATthree{\Linethree{1{+}x^3}{1{+}x}{1{+}x{+}x^2}}{\Linethree{x^2}{\!
1{+}x^2{+}x^3\!}{1{+}x^2}} {\Linethree{1{+}x{+}x^2\!}{1}{\!
1{+}x{+}x^2{+}x^3}}$ &
\doMATnine{
\MATnine{ \Linenine{1}{0}{0}{0}{0}{0}{0}{0}{0}}{
\Linenine{0}{1}{0}{1}{1}{1}{1}{1}{0}}{
\Linenine{0}{0}{1}{0}{0}{1}{0}{1}{0}}{
\Linenine{0}{0}{0}{1}{0}{0}{0}{1}{0}}{
\Linenine{0}{0}{0}{0}{1}{0}{1}{0}{0}}{
\Linenine{0}{0}{0}{0}{0}{1}{1}{0}{0}}{
\Linenine{0}{0}{0}{0}{0}{0}{1}{0}{0}}{
\Linenine{0}{0}{0}{0}{0}{0}{0}{1}{0}}{
\Linenine{0}{0}{0}{0}{0}{0}{0}{0}{1} }}{\MATnine{
\Linenine{0}{0}{0}{0}{0}{0}{0}{0}{0}}{
\Linenine{0}{0}{1}{1}{1}{1}{1}{1}{0}}{
\Linenine{0}{1}{1}{0}{0}{1}{0}{1}{0}}{
\Linenine{0}{1}{1}{0}{0}{1}{0}{1}{0}}{
\Linenine{0}{1}{0}{1}{1}{0}{1}{0}{0}}{
\Linenine{0}{1}{0}{1}{1}{0}{1}{0}{0}}{
\Linenine{0}{0}{1}{1}{1}{1}{1}{1}{0}}{
\Linenine{0}{1}{0}{1}{1}{0}{1}{0}{0}}{
\Linenine{0}{1}{1}{0}{0}{1}{0}{1}{0} }}
\\
$b_5$ &
$\MATthree{\Linethree{x^3}{1{+}x{+}x^2}{1}}{\Linethree{x{+}x^2}{x{+}x^2{+}x^3}{x}}{\Linethree{1}{x^2}{1{+}x^2{+}x^3}}$ &
\doMATnine{
\MATnine{ \Linenine{1}{0}{0}{0}{0}{0}{0}{0}{0}}{
\Linenine{0}{1}{0}{1}{1}{0}{0}{1}{0}}{
\Linenine{0}{0}{1}{0}{1}{0}{1}{1}{1}}{
\Linenine{0}{0}{0}{1}{0}{0}{0}{1}{0}}{
\Linenine{0}{0}{0}{0}{1}{0}{1}{0}{1}}{
\Linenine{0}{0}{0}{0}{0}{1}{1}{1}{1}}{
\Linenine{0}{0}{0}{0}{0}{0}{1}{0}{0}}{
\Linenine{0}{0}{0}{0}{0}{0}{0}{1}{0}}{
\Linenine{0}{0}{0}{0}{0}{0}{0}{0}{1} }}{\MATnine{
\Linenine{0}{0}{0}{0}{0}{0}{0}{0}{0}}{
\Linenine{0}{0}{1}{1}{1}{0}{0}{1}{0}}{
\Linenine{0}{1}{1}{0}{1}{0}{1}{1}{1}}{
\Linenine{0}{0}{1}{1}{1}{0}{0}{1}{0}}{
\Linenine{0}{1}{0}{1}{0}{0}{1}{0}{1}}{
\Linenine{0}{1}{1}{0}{1}{0}{1}{1}{1}}{
\Linenine{0}{0}{1}{1}{1}{0}{0}{1}{0}}{
\Linenine{0}{1}{1}{0}{1}{0}{1}{1}{1}}{
\Linenine{0}{0}{0}{0}{0}{0}{0}{0}{0} }}
\\
$b_6$ &
$\MATthree{\Linethree{x^2{+}x^3}{1}{x^2}}{\Linethree{1{+}x^2}{1{+}x^3}{1{+}x}}{\Linethree{x^2}{x{+}x^2}{x{+}x^2{+}x^3}}$ &
\doMATnine{
 \MATnine{ \Linenine{1}{0}{0}{0}{0}{0}{0}{0}{0}}{
\Linenine{0}{1}{0}{1}{0}{0}{1}{1}{1}}{
\Linenine{0}{0}{1}{1}{0}{1}{0}{0}{1}}{
\Linenine{0}{0}{0}{1}{0}{0}{0}{0}{1}}{
\Linenine{0}{0}{0}{0}{1}{0}{0}{0}{1}}{
\Linenine{0}{0}{0}{0}{0}{1}{1}{1}{1}}{
\Linenine{0}{0}{0}{0}{0}{0}{1}{0}{0}}{
\Linenine{0}{0}{0}{0}{0}{0}{0}{1}{0}}{
\Linenine{0}{0}{0}{0}{0}{0}{0}{0}{1} }}{\MATnine{
\Linenine{0}{0}{0}{0}{0}{0}{0}{0}{0}}{
\Linenine{0}{0}{1}{1}{0}{0}{1}{1}{1}}{
\Linenine{0}{1}{1}{1}{0}{1}{0}{0}{1}}{
\Linenine{0}{1}{1}{1}{0}{1}{0}{0}{1}}{
\Linenine{0}{1}{1}{1}{0}{1}{0}{0}{1}}{
\Linenine{0}{0}{1}{1}{0}{0}{1}{1}{1}}{
\Linenine{0}{1}{1}{1}{0}{1}{0}{0}{1}}{
\Linenine{0}{1}{0}{0}{0}{1}{1}{1}{0}}{
\Linenine{0}{1}{0}{0}{0}{1}{1}{1}{0} }}
\\ \hline
\end{tabular}
\end{center}
\caption{Generators of $\Gamma$ in the example of Section
\ref{sec:example}}\label{ExampleTab}
\end{figure}

\section{Glossary}\label{gloss}

$\AA(S)$ --- an $S$-algebra for any $R_T$-algebra $S$ with the
relations given in \Eq{rel} (see Section \ref{sec2})

$\B = \B_d(F)$ --- Bruhat-Tits building = $\PGL[d](F)/\PGL[d](\O)
= G(F)/K$

$\B^0$ --- the set of vertices of $\B$

$\B^i$ --- the $i$-skeleton of $\B$

$\B_I = \dc{\Gamma(I)}{G(F)}{K} \isom \Gamma/\Gamma(I)$

$d$ --- an integer $\geq 1$

$F = \F_q\db{y} = k_y$ (see Figure \ref{fi:fields_y})

$\tG'(S) = \mul{\AA(S)}$ (in particular $\tG'(R) = \mul{\AA(R)}$,
and likewise for the other groups defined over $R$ or other
$R_T$-algebras)

$G'(S) = \mul{\AA(S)}/\mul{S}$.

$\tG'_1(S)$ --- elements of $\tG'(S)$ of norm $1$

$G'_1(S)$ --- the image of $\tG'_1(S)$ under the map $\tG'(S) \ra
G'(S)$ (see Figure \ref{fi:GG})

$\tG(S) = \GL(S)$

$G(S) = \PGL(S)$

$\tG_1(S) = \SL(S)$

$G_1(S) = \PSL(S)$

$\con{\tG'}{S}{I}$ --- congruence subgroup of $\tG'(S)$; kernel of
$\tG'(S) \ra \tG'(S/I)$

$\con{G'}{S}{I}$ --- kernel of $G'(S) \ra G'(S/I)$ (for other
congruence subgroups, see Section \ref{sec:finq} and in particular
Figure \ref{fi:GG})

$\tilde{\Gamma} = \sg{b_u \suchthat u \in
\mul{\F_{q^d}}/\mul{\F_{q}}}$ (see Section \ref{sec:gamma} for
definition of $b_u \in \AA(R)$)

$\Gamma$ --- the image of $\tilde{\Gamma}$ under $\tG'(R) \ra
G'(R)$

$\Gamma(I) = \Gamma \cap \con{G'}{R}{I}$

$\Gamma_1 = \Gamma \cap G'_1(R)$ (see Figure \ref{GRsg})

$H = \con{G'}{R}{I}$

$H_1 = H \cap G'_1(R)$

$I = \ideal{p^s} \normali R$, $p \in \F_q[y]$ is a prime and $s
\geq 1$

$k = \F_q(y)$ (see Figure \ref{fi:fields_y})

$k_1$ = $\F_{q^d}(y)$ (see Figure \ref{fi:fields_y})

 $k_p$ --- the completion of $k$ with respect to $p$

$\bar{k} = \F_q(x)$

$L = R/I$, a finite local ring

$L_0 = R/J = R/\ideal{p}$, the residue field of $L$

$\bar{L} = \bar{R}/\bar{I}$

$\O = \F_q[[y]]$ --- the ring of integers of $F$

$\O_p$ --- the ring of integers of $k_p$

$q$ --- a fixed prime power, the order of the finite field $\F_q$

$R_T = \F_q[\ww,1/(1+\ww)] \sub \O$

$R =\F_q[\ww,1/\ww,1/(1+\ww)]$

$R_0 = \F_q[1/\ww]$

$\bar{R} = R[x]$

$\bar{R}_T = R_T[x]$

$\color \co \B^0 \ra \Z/d\Z$ --- the color function of the
building

$\nu_p$ --- valuation of $k$ for the prime polynomial $p \in
\F_q[y]$

$\yy{\nu}$ --- the minus degree valuation of $k$

$y$ ---   a transcendental variable over $\F_q$, generates $k$

\def\US{A.~Lubotzky, B.~Samuels and U.~Vishne}
\def\BIBPapI{\US, {\it Ramanujan complexes of type $\tilde{A_d}$},
Israel J. of Math., to appear.}
\def\BIBPapII{\US, {\it Explicit constructions of Ramanujan complexes},
European J. of Combinatorics, to appear.}
\def\BIBPapIII{\US, {\it Division algebras and non-commensurable isospectral manifolds},
preprint.}
\def\BIBPapIV{\US, {\it Isospectral Cayley graphs of some simple groups},
preprint.}

\end{document}